\documentclass{gtart_h}  

\def\ifplaintex{\expandafter\ifx\csname documentclass\endcsname\relax}

\ifplaintex 
\else
\fi

\expandafter\ifx\csname beginpicture\endcsname\relax
\expandafter\ifx\csname documentclass\endcsname\relax
\input pictex \else
\input prepictex \input pictex \input postpictex \fi\fi

\def\gt{{\mathsurround=0pt\it $\cal G\mskip-2mu$eometry \&\ 
$\cal T\!\!$opology}}        

\def\gtp{{\mathsurround=0pt\it $\cal G\mskip-2mu$eometry \&\ 
$\cal T\!\!$opology $\cal P\!$ublications}}  

\def\lognumber#1{\def\thelognumber{#1}}
\def\volumenumber#1{\def\thevolumenumber{#1}}
\def\papernumber#1{\def\thepapernumber{#1}}
\def\volumeyear#1{\def\thevolumeyear{#1}}

\def\pagenumbers#1#2{\def\startpage{#1}\def\finishpage{#2}}
\def\published#1{\def\publishdate{#1}}
\def\proposed#1{\def\theproposer{#1}}
\def\seconded#1{\def\theseconders{#1}}
\def\received#1{\def\receiveddate{#1}}

\def\accepted#1{\def\accepteddate{#1}}

\long\def\asciiabstract#1{\long\def\theasciiabstract{#1}}

\let\\\par\let\thelognumber\relax
\let\thevolumenumber\relax\let\thepapernumber\relax
\let\thevolumeyear\relax\let\thesamplenumber\relax\let\startpage\relax
\let\finishpage\relax\let\publishdate\relax\let\receiveddate\relax
\let\reviseddate\relax\let\accepteddate\relax\let\theasciititle\relax
\let\theasciiauthors\relax
\let\theasciiabstract\relax
\let\theasciiemail\relax\let\theshortauthors\relax\let\theshorttitle\relax

\long\def\maketitlep{   

\count0=\startpage

\gt\hfill      
\beginpicture
\setcoordinatesystem units <0.33truein, 0.33truein> point at 2.2 0.9
\setplotsymbol ({$\cal G$})
\plotsymbolspacing=9truept
\circulararc 315 degrees from 0 1 center at 0 0
\setplotsymbol ({$\cal T$})
\circulararc 315 degrees from 1 -1 center at 1 0
\endpicture
\break
{\small\ifx\thesamplenumber\relax 
Volume \else Sample
\fi\thevolumenumber\ (\thevolumeyear)
\startpage--\finishpage\nl
Published: \publishdate}
\vglue 0.5truein plus 0.4fil minus 0.1truein

{\parskip=0pt\leftskip 0pt plus 1fil\def\\{\par\smallskip}{\ifplaintex\large
\else\Large\fi\bf\thetitle}\par\medskip}   

\vglue 0pt plus 0.1fil 

{\parskip=0pt\leftskip 0pt plus 1fil\def\\{\par}{\sc\theauthors}
\par\medskip}

\vglue 0pt plus 0.1fil 

{\small\parskip=0pt\let\newline\\
{\leftskip 0pt plus 1fil\def\\{\par}{\sl\theaddress}\par}
\expandafter\ifx\theemail\relax    
\relax\else\vglue 5pt plus 0.02fil minus 2pt\def\\{\stdspace{\rm 
and}\stdspace} 
\cl{Email:\stdspace\tt\theemail}\fi
\ifx\theurl\relax                  
\relax\else\vglue 5pt plus 0.02fil minus 2pt\def\\{\stdspace{\rm 
and}\stdspace}
\cl{URL:\stdspace\tt\theurl}\fi\par}

\vglue 7pt plus 0.3fil minus 3pt

{\bf Abstract}
\vglue 5pt plus 0.1fil minus 2pt

\theabstract

\vglue 7pt plus 0.3fil minus 3pt

{\bf AMS Classification numbers}\quad Primary:\quad \theprimaryclass

Secondary:\quad \thesecondaryclass

\vglue 5pt plus 0.3fil minus 2pt

{\bf Keywords:}\quad \thekeywords

\vglue 10pt plus 0.5fil minus 5pt

{\small  Proposed: \theproposer\hfill Received: \receiveddate\nl
Seconded: \theseconders\hfill 
\ifx\reviseddate\relax                         
Accepted: \accepteddate                        
\else
Revised: \reviseddate                          
\fi}
\eject
}       

\let\maketitlepage\maketitlep
\let\maketitle\maketitlepage

\font\phead=cmsl9 scaled 950
\font=cmsl9 scaled 1050
\font\pnum=cmbx10 scaled 913
\font\lnum=cmbx10 
\font\pfoot=cmsl9 scaled 950
\font=cmsl9 scaled 1050
\ifplaintex
\headline{\vbox to 0pt{\vskip -4.5mm\line{\small\phead\ifnum
\count0=\startpage ISSN 1364-0380 (on line)
1465-3060 (printed) \hfill {\pnum\folio}\else\ifodd\count0\def\\{ }%
\ifx\theshorttitle\relax\thetitle\else\theshorttitle\fi\hfill{\pnum\folio}
\else\def\\{ and }{\pnum\folio}\hfill\ifx\theshortauthors\relax\theauthors
\else\theshortauthors\fi\fi\fi}\vss}}
\footline{\vbox to 0pt{\vglue 0mm\line{\small\pfoot\ifnum\count0=\startpage
\copyright\ \gtp\hfill\else
\gt, Volume \thevolumenumber\ (\thevolumeyear)\hfill\fi}\vss
}}
\else
\makeatletter
\def\@oddhead{{\small\lhead\ifnum\count0=\startpage ISSN 1364-0380 (on line)
1465-3060 (printed) \hfill {\lnum\number\count0}\else\ifodd\count0
\def\\{ }\ifx\theshorttitle\relax \thetitle \else\theshorttitle\fi\hfill
{\lnum\number\count0}\else\def\\{ and }{\lnum\number\count0}
\hfill\ifx\theshortauthors\relax 
\theauthors\else\theshortauthors\fi\fi\fi}}\def\@evenhead{\@oddhead}
\def\@oddfoot{\small\lfoot\ifnum\count0=\startpage\copyright\ \gtp\hfill\else
\gt, Volume \thevolumenumber\ (\thevolumeyear)\hfill\fi}
\def\@evenfoot{\@oddfoot}
\makeatother
\fi

\newwrite\gtoutfile
\long\gdef\makeheadfile{  
{\def\\{, }\def\s{ }
\immediate\openout\gtoutfile head.xxx
\immediate\write\gtoutfile{Proxy-for: \ifx\theasciiauthors\relax
\theauthors\else\theasciiauthors\fi\s<\ifx\theasciiemail\relax\theemail\else\theasciiemail\fi>}
\immediate\write\gtoutfile{\noexpand\\}
\immediate\write\gtoutfile{Authors: \ifx\theasciiauthors\relax
\theauthors\else\theasciiauthors\fi}
{\def\\{ }\immediate\write\gtoutfile{Title: \ifx\theasciititle\relax
\thetitle\else\theasciititle\fi}}
\immediate\write\gtoutfile{Subj-class: GT or SG or MG etc}
\immediate\write\gtoutfile{MSC-class: \theprimaryclass\ifx\thesecondaryclass\relax\else, \thesecondaryclass\fi}
\immediate\write\gtoutfile{Journal-ref: Geom. Topol. \thevolumenumber
(\thevolumeyear) \startpage-\finishpage}
\immediate\write\gtoutfile{Comments: Published by Geometry and Topology at}
\immediate\write\gtoutfile{\s\s http://www.maths.warwick.ac.uk/gt/GTVol\thevolumenumber/paper\thepapernumber.abs.html}
\immediate\write\gtoutfile{\noexpand\\}
\immediate\write\gtoutfile{}
\ifx\theasciiabstract\relax
\immediate\write\gtoutfile{\theabstract}\else
\immediate\write\gtoutfile{\theasciiabstract}\fi
\immediate\write\gtoutfile{}
\immediate\write\gtoutfile{\noexpand\\}
\immediate\write\gtoutfile{}
\immediate\closeout\gtoutfile}}  

\def\maketitlepage{\maketitlep\makeheadfile}
\let\maketitle\maketitlepage

\lognumber{577}
\received{12 March 2005}
\volumenumber{9}\papernumber{55}\volumeyear{2005}
\pagenumbers{2359}{2394}   
\accepted{16 December 2005}
\proposed{Benson Farb}
\seconded{Jean-Pierre Otal, Steven Ferry}
\published{26 December 2005}

\usepackage{graphicx,psfrag,amsmath,amssymb}

\def\psfraga <#1,#2> #3#4{%
\psfrag {#3}{\smash{\rlap{\kern #1 \raise #2\hbox{#4}}}}}

\def\figref#1{\hyperlink{#1anchor}{Figure~\ref*{#1}}}
\def\anchor#1{\noindent\hypertarget{#1anchor}{\smash{$\phantom{99}$}}\newline}

\newtheorem{theorem}{Theorem}

\newtheorem{conjecture}[theorem]{Conjecture}
\newtheorem{corollary}[theorem]{Corollary}

\newtheorem{lemma}[theorem]{Lemma}

\newtheorem{proposition}[theorem]{Proposition}

\begin{document}
\title{On the dynamics of isometries}
\author{Anders Karlsson}
\address{Mathematics Department, Royal Institute of Technology\\100 44 
Stockholm, Sweden}
\email{akarl@math.kth.se}

\begin{abstract}
We provide an analysis of the dynamics of isometries and semicontractions
of metric
spaces. Certain subsets of the boundary at infinity play a
fundamental role and are identified completely for the standard
boundaries of CAT(0)--spaces, Gromov hyperbolic spaces, Hilbert
geometries, certain pseudoconvex domains, and partially for
Thurston's boundary of Teichm\"{u}ller spaces.
We present several rather general results concerning groups
of isometries, as well as the proof of other more specific
new theorems,
for example concerning the existence of free
nonabelian subgroups in CAT(0)--geometry, iteration of holomorphic maps,
a metric Furstenberg lemma, random walks on groups, noncompactness of
automorphism groups of convex cones,
and boundary behaviour of Kobayashi's metric.
\end{abstract}
\asciiabstract{%
We provide an analysis of the dynamics of isometries and
semicontractions of metric spaces. Certain subsets of the boundary at
infinity play a fundamental role and are identified completely for the
standard boundaries of CAT(0)-spaces, Gromov hyperbolic spaces,
Hilbert geometries, certain pseudoconvex domains, and partially for
Thurston's boundary of Teichmueller spaces.  We present several
rather general results concerning groups of isometries, as well as the
proof of other more specific new theorems, for example concerning the
existence of free nonabelian subgroups in CAT(0)-geometry, iteration
of holomorphic maps, a metric Furstenberg lemma, random walks on
groups, noncompactness of automorphism groups of convex cones, and
boundary behaviour of Kobayashi's metric.}

\primaryclass{37B05, 53C24, }
\secondaryclass{22F50, 32H50}
\keywords{Metric spaces, isometries, nonpositive curvature, 
Kobayashi metric, random walk}
\maketitlepage

\section{Introduction}

The notion of a metric space was introduced in 1906 by M Fr\'echet.
Although a systematic study of metric spaces from the point of
view of point-set topology was subsequently undertaken,
the most natural morphisms, isometries
and semicontractions, seem to have received less attention.
To be sure, there are special topics that have inspired deep
investigations: Euclidean and hyperbolic geometry, extensions of
the contraction mapping principle, iteration of holomorphic maps
as well as, in more recent years, CAT(0)--spaces and Gromov
hyperbolic groups. But in contrast to the category of topological
vector spaces and continuous linear operators, a basic general
text on metric spaces and semicontractions seems to be absent.
Note that there are a number of contexts in for example geometry,
topology, complex analysis in one and several variables, Lie
theory, ergodic theory and group theory, where metrics and
semicontractions arise. Some of these will be recalled in more 
detail later on.

This paper presents a general and unified theory of the dynamics
of semicontractions and (groups of) isometries. It studies and
exploits (generalized) halfspaces and their limits, \emph{the
stars at infinity}. These subsets are of fundamental importance
for the dynamics of isometries and provide moreover a convenient
framework for asymptotical geometric information, and should
therefore be of interest to the subjects of Riemannian and metric
geometry. Even though halfspaces are classical in the definition
of Dirichlet fundamental domains and appear particularly in the
literature on Kleinian groups, it seems they have not been
systematically considered previously. The stars relate well to
standard concepts such as Tits geometry of CAT(0)--spaces,
Thurston's boundary of Teichm\"{u}ller space, hyperbolicity of
metric spaces, strict pseudoconvexity, the face lattice of convex
domains, rank 1 isometries, etc.

In the theory of word hyperbolic groups, the study of how the
group acts on its boundary plays an important role.
Our generalizations of hyperbolic phenomena bringing in the stars
and their
incidence geometry, are perhaps also interesting in light of
Mostow's proof of strong rigidity in the higher rank case.

Several of the results obtained are new even in areas which have
been much studied, for example CAT(0)--geometry or boundary
behaviour in complex domains and holomorphic maps.  Let us
highlight a few of these results:

\begin{theorem}
Let $X$ be a proper\/ {\rm CAT(0)}--space. Assume $g_n$ is a sequence of
isometries
such that $g_{n}x_{0}\rightarrow\xi^{+}\in\partial X$ and $g_{n}^{-1}%
x_{0}\rightarrow\xi^{-}\in\partial X$. Then for any
$\eta\in\overline{X}$ with $\angle(\eta,\xi^{-})>\pi/2$ we have
that
$$
g_{n}\eta\rightarrow\{\zeta:\angle(\xi^{+},\zeta)\leq\pi/2\}
$$
and the convergence is uniform outside neighborhoods of
$S(\xi^{-})$.
\end{theorem}

This is completely new in that it equally well deals with
parabolic isometries. Previously, mainly iterates of a single
hyperbolic isometry could be treated, in which case a lemma of
Schroeder \cite{BGS 85} generalized by Ruane \cite{Rua01} to
include also singular CAT(0)--spaces actually gives more
information. For several other new corollaries on groups acting on
CAT(0)--spaces, see section \ref{seccat0}. Next, the following
describes a novel phenomenon for simple random walks on \emph{any}
finitely generated nonamenable group:

\begin{theorem}
Let $\Gamma$ be a nonamenable group generated by a finite set $S$
and consider the random walk defined by the uniform distribution
on $S\cup S^{-1}$. For almost every trajectory there is a time
after which every finite collection of halfspaces defined by the
trajectory intersect nontrivially.
\end{theorem}

For more discussion and explanations, see subsection \ref{ssecrw}.
Every holomorphic map is in a sense a semicontraction and taking
advantage of this we will obtain the following new Wolff--Denjoy
theorem:

\begin{theorem}
Let $X$ be a bounded $C^2$--domain in $\mathbb{C}^n$ which is
complete in the Kobayashi metric satisfying the boundary estimate
(\ref{est1}) in subsection \ref{ssecmetdist}. Let $f:X\rightarrow
X$ be a holomorphic map. Then either the orbit of $f$ stays away
from the boundary or there is a unique boundary point $\xi$ such
that
$$
\lim_{m\rightarrow\infty} f^m(z)=\xi
$$
for any $z\in X$.
\end{theorem}

Examples include real analytic pseudoconvex domains in which case
the theorem for $n=2$ was proved by Zhang and Ren in \cite{ZR95}.
Finally, we mention:

\begin{theorem}
Any polyhedral cone with noncompact automorphism group has
simplicial diameter at most 3.
\end{theorem}

Here one should note that in dimension $2$ a rather complete
result concerning which convex sets have infinite automorphism
group can be found in de la Harpe's paper \cite{dlH93}.

I would like to thank Bruno Colbois and Alain Valette
for inviting me to spend a very pleasant and productive year at
Universit\'e de Neuch\^atel during which a large part of this work
was written. Support from the Swiss National Science Foundation
grant 20-65060.01 and the Swedish Research Council grant
2002-4771 are gratefully acknowledged.

\section*{Part I\qua General theory}
\addcontentsline{toc}{section}{Part I: General Theory}

\section{Halfspaces and stars at infinity}

\subsection{Definitions} 
\begin{figure}[ht!]\small\anchor{fig:parab}
\psfrag{x}{$x_0$}
\psfraga <-3pt,0pt> {W}{$W$}
\psfraga <-3pt,0pt> {H}{$H(W)$}
\cl{\includegraphics{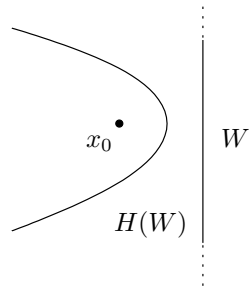}} \caption{The halfspace
defined by a line $W$ in $\mathbb{R}^2$ is the region containing
$W$ bounded by a parabola.}\label{fig:parab}
\end{figure}

Let $X$ be a metric space. For a
subset $W$ of $X$ we let
\[
d(x,W)=\inf_{w\in W}d(x,w).
\]
Fix a base point $x_{0}$. We define the \emph{halfspace} defined
by the subset
$W$ and the real number $C$ to be%
\[
H(W,C)=H^{x_{0}}(W,C):=\{z:d(z,W)\leq d(z,x_{0})+C\}.
\]
We use the notation $H(W):=H(W,0)$ and for two points $x$ and $y$
in $X$ we let $H_{y}^{x}=\{z:d(z,y)\leq d(z,x)\},$ so
$H_{y}^{x_{0}}=H(\{y\},0)$. Note that the latter sets define
halfspaces in the more standard sense when $X$ is a Euclidean or
real hyperbolic space.  See \figref{fig:parab}.

Let $X$ be a complete metric space. By a \emph{bordification} of
$X$ we here mean a Hausdorff topological space $\overline{X}$ with $X$
embedded as an open dense subset. The \emph{boundary }is $\partial
X=\overline{X}\setminus X$. If $\overline X$ is compact we refer
to it as a \emph{compactification}. We define $d(x,\xi)=\infty$
for any $x\in X$ and $\xi\in\partial X$ (which is consistent with
the completeness of $X$) and extend the definition of $d(x,W)$ for
$W\subset\overline{X}$ in the expected way. 

A metric space is
\emph{proper} if every closed ball is compact. Recall that every proper
metric space $X$ has a (typically nontrivial) metrizable
Isom($X$)--compactification $\overline{X}^h$
by horofunctions: $X$ is embedded into the space of continuous functions
$C(X)$ with the topology of uniform convergence on bounded sets via
$$
x\mapsto d(x,\cdot)-d(x,x_0)
$$
The closure of the image now defines the compactification,
see \cite{BGS 85}, \cite{Bal95}, and \cite{BH99} for more details.

\medskip
\textbf{Example}\qua
Another general compactification, \emph{the end compactification}, was introduced by
Freudenthal. Here let $X$
be path connected, proper metric space and define the following equivalence relation on
the set of proper rays from $x_0$. Two rays are equivalent if for any
compact set $K$ in $X$, the two rays are eventually contained in the 
same path connected components of $X\setminus K$. The equivalence classes
of proper rays union $X$ with the natural topologization constitute
the compactification of $X$.
See \cite{BH99} for more details.
\medskip

Let $\mathcal{V}_{\xi}$ denote the collection of open
neighborhoods in $\overline{X}$ of a boundary point $\xi$. The
\emph{star based at $x_0$} of a point $\xi\in\partial X$ is
$$
S^{x_{0}}(\xi):=\bigcap_{V\in\mathcal{V_\xi}}\overline{H(V)},
$$
where the closures are taken in $\overline{X}$, and the
\emph{star} of $\xi$ is
$$
S(\xi):=\overline{ {\displaystyle\bigcup\limits_{C\geq0}}
\bigcap_{V\in\mathcal{V_\xi}}\overline{H(V,C)}}.
$$
The latter definition in particular removes an a priori dependence
of $x_0$ as will be clear later on. Note also that because of the
monotonicity built into the definition of $H,$ we may restrict
$\mathcal{V}_{\xi}$ to some fundamental system of neighborhoods of
$\xi$.

\begin{figure}[ht!]\anchor{fig:defstar2}\small
\psfraga <-4pt,0pt> {H}{$H(V,C)$}
\psfrag {x}{$x_0$}
\psfraga <-4pt,0pt> {xi}{$\xi$}
\psfraga <-4pt,0pt> {V}{$V$}
\psfraga <-4pt,0pt> {S}{$S(\xi)$}
\cl{\includegraphics{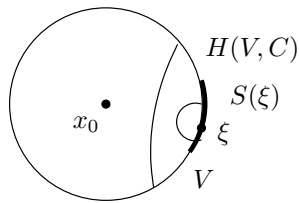}}
\caption{The definition of $S(\xi)$}\label{fig:defstar2}
\end{figure}

We introduce the \emph{star-distance}: Let $s$ be the largest
metric on $\partial X$ taking values in $[0,\infty]$ such that
$s(\xi,\eta)=0$ if $S(\xi)=S(\eta)$, and $s(\xi,\eta)=1$ if at
least one of $\xi\in S(\eta)$ or $\eta\in S(\xi)$ holds. More
explicitly, $s(\xi,\eta)$ equals the minimum number $k$ such that
there are points $\gamma_i$ with $\gamma_0=\xi$, $\gamma_k=\eta$,
and $s(\gamma_i,\gamma_{i+1})=1$ for all $i$.  See \figref{fig:defstar2}.

\medskip
\textbf{Example}\qua
Let $X$ be the Euclidean space $\mathbb{R}^n$ and $\partial X$
the visual sphere at infinity being the space of geodesic rays from the 
origin. Then $S(\xi)$ as well as $S^{x_0}(\xi)$ is the hemishpere 
in $\partial X$ centered at $\xi$. For a generalization of this
example, see Proposition \ref{propcat0}. Hence the visual sphere has
stardiameter 2. Indeed, for all points $\eta \in S(\xi)$ different from
$\xi$ we have $s(\eta,\xi)=1$ and for the points $\zeta$ 
outside this star $s(\zeta,\xi)=2$.

\medskip
\textbf{Example}\qua
Let $X$ be a proper and path connected metric space and $\partial X$
the space of ends as defined above. Consider two nonequivalent proper
rays denoted $\eta$ and $\xi$ (with an abuse of notation). For 
any two small disjoint neighborhoods of these points, there is
a compact set $K$ which separate these two neighborhoods in the
sense that they are in different path components of
$X\setminus K$. Therefore any path between the two neighborhoods
must pass through $K$, and it follows that $\eta$ is not in $S(\xi)$.
Hence $S(\xi)=S^{x_0}(\xi)=\{\xi\}$ and $s(\xi,\eta)=\infty$
for any two distinct boundary points. See section \ref{sechyp} for
further examples with this kind of ``trivial" or ``hyperbolic"
star geometry.
\smallskip

It does not seem clear whether, or when, $\xi\in S(\eta)$ implies
$\eta\in S(\xi)$. Let
$$S^{\vee}(\xi)=\{\eta :\xi\in S(\eta)\},$$
and we say that the bordification is \emph{star-reflexive} when
$S(\xi)=S^{\vee}(\xi)$ for all $\xi$. The examples below turn out
to have this property.

The \emph{face} of a subset $A$ of $\partial X$ is the
intersection of all stars containing $A$. The face of the empty
set is defined to be the empty set.
We define for a subset $A\subset \partial X$ the sets
$$
S(A)=\bigcap_{a\in A} S(a).
$$
and similarily for $S^{\vee}(A)$.

By the notation
$x_{n}\rightarrow S,$ where $x_{n}$ is a sequence of points and
$S$ a set, we mean that for any neighborhood $U$ of $S$ we have
$x_{n}\in U$ for all sufficiently large $n$.

\subsection{Some lemmas}
\begin{lemma}
\label{lemS}For any $\xi\in\partial X$, the sets $\overline{H(V)}$ for
$V\in\mathcal{V}_{\xi}$ \ contain $V$ and $\xi\in
S^{x_{0}}(\xi)\subset S(\xi)\subset\partial X$. If $\partial X$ is
compact, then for every neighborhood $U$ of $S^{x_{0}}(\xi)$ there
is a neigborhood $V$ of $\xi$ such that $\overline{H(V)}\subset
U$.
\end{lemma}

\begin{proof}
Note that $V\subset\overline{H(V)}$. Indeed, first observe that
$V\cap X\subset H(V)$ because $d(v,V)=0$ for any $v\in V$.
Secondly, note that for any $v\in V$ and any open neighborhood $U$
of $v,$ $U\cap V$ is again an open neighborhood and every open set
in $\overline{X}$ has to intersect $X$. Finally, $S^{x_{0}}(\xi)$
is nonempty because $\xi$ is contained in every $V,$ and
$S^{x_{0}}(\xi)\subset S(\xi)\subset\partial X$ since $d(V,x_{0})$
is unbounded for $V\in\mathcal{V}_{\xi}$.

Let $U$ be a neighborhood of $S^{x_0}(\xi)$. We may assume that $U$
is open and so $U^{c}$ is compact. Consider a fundamental system of
neighborhoods of $\xi$. Suppose for any $V$ in this system it holds
that $\overline{H(V)}\cap U^{c}\neq\emptyset$. Becuase of the 
monotonicity of halfspaces we hence have a decreasing, nested system
of closed sets $\overline{H(V)}\cap U^{c}$ inside $U^{c}$. By
compactness we get ${\displaystyle\bigcap} \overline{H(V)}\cap U^{c}
\neq\emptyset)$. This is a contradiction to $S^{x_0}(\xi)\subset U$, 
and proves the last assertion of the lemma.
\end{proof}

Note that if $z_{n}\rightarrow \xi$ and $d(z_{n},y_{n})<C$ then
every limit point of $y_{n}$ belongs to $S^{x_{0}}(\xi)$. A
priori, $S^{x_{0}}(\xi)$ depends on $x_{0}$ although in the
examples below this turns out not to be the case. On the other hand:

\begin{lemma}
\label{lemSindep}The sets $S(\xi)$ are independent of the base point $x_{0}$.
If $z_{n}\rightarrow\xi\in\partial X$, $d(z_{n},y_{n})<C$ and
$y_{n}\rightarrow\eta$, then $S(\xi)=S(\eta)$. Moreover, $\xi$ and
$\eta$ belong to the same stars.
\end{lemma}

\begin{proof} The first statement follows from
$$
H^{x_{0}}(W,C-d(x,x_{0}))\subset H^{x}(W,C)\subset
H^{x_{0}}(W,C+d(x,x_{0})),
$$
and because of the increasing union over $C\geq0$ in the
definition of $S(\xi)$. The other two claims hold for similar
reasons.
\end{proof}

\begin{lemma}
\label{lemscont}Assume that $\overline X$ is sequentially compact and
that $S(\xi)=S^{x_{0}}(\xi)$ for every $\xi\in\partial X$. Let
$\xi_{n}$ and $\eta_{n}$ be two sequences in $\partial X$
converging to $\xi$ and $\eta$, respectively. If
$s(\xi_{n},\eta_{n})>0$ for all $n$, then
$$
s(\xi,\eta)\leq\liminf_{n\rightarrow\infty}s(\xi_{n},\eta_{n}).
$$
\end{lemma}

\begin{proof} By the assumption we can work with the
$S^{x_{0}}$--stars. It is enough to consider
$s(\xi_{n},\eta_{n})=1$ for all $n$, because of the sequential
compactness and the way $s$ is defined. Moreover, we may suppose
that $\xi_{n}\in S(\eta_{n})$ for all $n$. Hence
$\xi_{n}\in\overline{H(V)}$ for every neighborhood $V$ of
$\eta_{n}$. Given a neighborhood $U$ of $\eta$, there is a $N$
such that $U$ is also a neighborhood of $\eta_{n}$ for $n\geq N$.
We therefore have that $\xi_{n}\in\overline{H(U)}$ for all $n\geq
N$, and hence also $\xi\in\overline{H(U)}$. Because $U$ was
arbitrary, we have that $\xi\in S(\eta)$ and so $s(\xi,\eta)\leq1$
as required.
\end{proof}

\section{Dynamics of isometries}

\subsection{Definitions} Let $X$ be a metric space. A
\emph{semicontraction} $f$ is a map $f:X\rightarrow X$ such
that
$$
d(f(x),f(y))\leq d(x,y)
$$
for every $x$, $y\in X$.
An \emph{isometry} is here an isomorphism in this categroy, which
means it is a distance preserving bijection.

A subset $D$ of semicontractions is called \emph{bounded} (resp.\
\emph{unbounded}) if $Dx_{0}$ is a bounded (resp.\ an unbounded) set.
A single semicontraction $f$ is called \emph{bounded} (resp.\
\emph{unbounded}) if $\{f^{n}\}_{n>0}$ is bounded (resp.\
unbounded). Note that these definitions are independent of
$x_{0}$.

If the action of the isometries of $X$ extends to an action by
homeomorphisms of $\overline{X}$ we call the bordification an
\emph{Isom}($X$)\emph{--bordification}. Note that when $X$ is a
proper metric space, the horofunction compactification
$\overline{X}^h$ is a (almost always nontrivial) metrizable
Isom($X$)--compactification (see the previous section).

Under the assumption that $\overline{X}$ is an
Isom($X$)--bordification, the
isometries of $X$ act on the stars $S(\xi)$ as can be seen from:%
\begin{eqnarray*}
gH(W,C)&=&\{z:d(g^{-1}z,W)\leq d(g^{-1}z,x_{0})+C\}\\
&=&\{z:d(z,gW)\leq d(z,gx_{0})+C\},
\end{eqnarray*}
which is included in $H(gW,C+d(x_{0},gx_{0}))$ and contains $H(gW,C-d(x_{0}%
,gx_{0})).$ Hence we have $gS(\xi)=S(g\xi)$ and it is plain that
$g$ preserves star distances. Note that we also have an action on
the faces.

\subsection{A contraction lemma\label{sseccontr}} The
following observation lies behind the construction of Dirichlet
fundamental domains (see eg \cite{Rat94}): For any isometry $g$
it holds that
$$
g(H_{x}^{g^{-1}y})=H_{gx}^{y}.
$$
This leads to a contraction lemma, which in spite of its
simplicity and fundamental nature, we have not been able to locate
in the literature:

\begin{lemma}
\label{lemcontract}Let $g_{n}$ be a sequence of isometries such that
$g_{n}x_{0}\rightarrow\xi^{+}$ and
$g_{n}^{-1}x_{0}\rightarrow\xi^{-}$ in a bordification
$\overline{X}$ of $X$. Then for any neighborhoods $V^{+}$ and
$V^{-}$ of $\xi^{+}$ and $\xi^{-}$ respectively, there exists
$N>0$ such that
$$
g_{n}(X\setminus H(V^{-}))\subset H(V^{+})
$$
for all $n\geq N$.
\end{lemma}

\begin{proof} Given neighborhoods $V^{+}$ and $V^{-}$ as in the
statement, by assumption there is an $N$ such that $g_{n}x_{0}\in
V^{+}$ and $g_{n}^{-1}x_{0}\in V^{-}$ for every $n\geq N$. For any
$z\in X$ outside $H(V^{-})$, so $d(z,v)>d(z,x_0)$ for every $v\in
V^{-}$, we have
$$
d(g_{n}z,V^{+})\leq d(g_{n}z,g_{n}x_{0})=d(z,x_{0})
<d(z,g_{n}^{-1}x_{0})=d(g_{n}z,x_{0})
$$
for every $n\geq N$.
\end{proof}

Here is a version of the contraction phenomenon when the
isometries act on the boundary:

\begin{proposition}\label{propcontrS}
\label{propcontr}Assume that $\overline{X}$ is an
Isom($X$)--compactification. Let $g_{n}$ be a sequence of
isometries such that $g_{n}x_{0}\rightarrow \xi^{+}$ and
$g_{n}^{-1}x_{0}\rightarrow\xi^{-}$ in $\overline{X}.$ Then for
any $z\in\overline{X}\setminus S^{x_{0}}(\xi^{-})$,
$$
g_{n}z\rightarrow S^{x_0}(\xi^{+}).
$$
Moreover, the convergence is uniform outside neighborhoods of
$S^{x_0}(\xi^{-})$.
\end{proposition}

\begin{proof}
Since $z$ does not belong to $S^{x_{0}}(\xi^{-})$
there is some neighborhood $V^{-}$ of $\xi^{-}$ such that
$z\notin\overline{H(V^{-})}$. As the latter is a closed set, there
is an open neighborhood $U$ of $z$ disjoint from
$\overline{H(V^{-})}$. Given a neighborhood $V^{+}$ of $\xi^{+}$
we therefore have for all sufficiently large $n$ that $g_{n}(U\cap
X)\subset H(V^{+})$ for all $n>N$. Since $g_{n}$ are
homeomorphisms we have that $g_{n}z\subset \overline{H(V^{+})}$ as
required. The conclusion now follows in view of Lemma \ref{lemS}.
\end{proof}

In some cases, for example if $z\in S(\eta)$ for some $\eta\notin
S(\xi^-)$, one can say more in view that
the isometries preserve stars. 

\subsection{Individual semicontractions} Let $f$ be a
semicontraction of a complete metric space $X$ and let
$\overline{X}$ be a bordification of $X$.
The \emph{limit set} of the $f$--orbit of $x_0$ is
$$L^{x_0}(f):=\overline{\{f^{n}(x_{0})\}}_{n>0}\cap\partial X,$$
which necessarily is empty if $f$ is bounded.

\begin{proposition}\label{propfixed} Let $g$ be an (unbounded)
isometry and $\overline{X}$ an Isom($X$)--bordification. Then $g$
fixes the star and the face of every point $\xi\in L^{x_0}(g)$, 
that is, $S(g\xi)=S(\xi)$ and $F(g\xi)=F(\xi)$.
Moreover, the subset $F(g)$ defined below, is also fixed by $g$.
\end{proposition}

\begin{proof} Since by continuity
$$
g\xi=g(\lim_{k\rightarrow\infty}g^{n_{k}}x_{0})=\lim_{k\rightarrow\infty
}g^{n_{k}}(gx_{0})
$$
we have that $S(g\xi)=S(\xi)$ in view of Lemma \ref{lemSindep}. If
$\xi\in S(\eta),$ then $g\xi\in S(g\eta)$ and again we have
$\xi\in S(g\eta).$ Since $g$ is a bijection, the final part of the
proposition follows.
\end{proof}

Let $a_n=d(f^n(x_0),x_0).$ A subsequence $n_i\rightarrow\infty$ is
called \emph{special} for
$f$ if $a_{n_i}\rightarrow\infty$ and there is a constant $C\geq 0$ 
such that $a_{n_{i}}>a_{m}-C$
for all $i$ and $m<n_{i}$. Note that being special clearly passes
to subsequences and by the triangle inequality it is independent
of $x_0$ (see (\ref{triangle}) below). Moreover, special subsequences
are invariant under the
shift $\{n_i\}\mapsto\{n_i+N\}$, where $N$ is some fixed integer.

Let $A^{x_0}(f)$ denote the limit points of $f^n(x_0)$ along the
special subsequences. The \emph{characteristic set} $F(f)$ of $f$
is the face of $A^{x_0}(f)$. (It may of course happen that
$A^{x_0}(f)=\emptyset$, in which case $F(f):=\emptyset$.)

\begin{theorem}\label{thmindivid} Assume that $X$ is proper and that
$\overline{X}$ is a sequentially compact bordification of $X$.
To any semicontraction $f$, the subset $F(f)\subset\partial X$ is
canonically associated to $f$. It holds that
$F(f)=\emptyset$ if and only if $f$ is bounded, and that if
$F(f)\neq\emptyset$, then every $f$--orbit
accumulates only at $\partial X$. Moreover, for any $x_0\in X$
$$
L^{x_0}(f)\subset S^{\vee}(F(f)). 
$$
If in addition $\overline{X}$ is star-reflexive, then
$$
L^{x_0}(f)\subset S(F(f)).
$$
\end{theorem}

\begin{proof}
 From the triangle inequality we get
\begin{equation}\label{triangle}
|d(g^{k}x,x)-d(g^{k}x_{0},x_{0})|\leq2d(x,x_{0}),
\end{equation}
which implies in view of Lemma \ref{lemSindep} that $F(f)$ is
independent of $x_{0}$. By completeness, if $f$ is bounded, then
$F(f)=\emptyset$. The converse is proved below. Calka's theorem
\cite{Cal84} asserts that if there is a bounded subsequence
of the orbit, then in fact the whole orbit is bounded.

Now suppose $f$ is unbounded and let $n_i$ be a special
sequence for $f$ (it is obvious, see \cite{Ka01}, that
special subsequences exist if and only if $f$ is unbounded) and
such that $f^{n_{i}}(x_{0})$ converges to some point $\xi\in\partial
X$. Observe that for any positive $k<n_{i}$ it holds that
$$
d(f^{n_{i}}(x_{0}),f^{k}(x_{0}))\leq d(f^{n_{i}-k}(x_{0}),x_{0})
=a_{n_{i}-k}<a_{n_{i}}+C=d(f^{n_{i}}(x_{0}),x_{0})+C.
$$
Now suppose we have a convergent sequence
$f^{k_{j}}x_{0}\rightarrow\eta \in\partial X$, which means that
given a neighborhood $V$ of $\eta$, we can find $j$ large so that
$f^{k_{j}}x_{0}\in V$. Now from the above inequality we get that
for all large enough $i$
$$
f^{n_{i}}x_{0}\in H(\{f^{k_{j}}x_{0}\},C)\subset H(V,C).
$$
Therefore $\xi\in\overline{H(V,C)}$ and since $V$ was an arbitrary
neighborhood we have $\xi\in S(\eta)$. (Note that in particular
this means that $A^{x_0}(f)$ and $F(f)$ are nonempty.) Finally,
assuming star-reflexivity we have showed that $\eta\in S(\xi)$ for
every special limit point $\xi$.
\end{proof}

\begin{figure}[ht!]\small\anchor{fig:orbit}
\psfraga <-5pt,0pt> {F}{$F(f)$}
\psfraga <-2pt,0pt> {x}{$x_0$}
\cl{\includegraphics{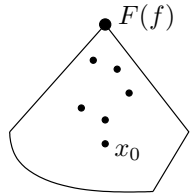}}
\caption{The orbit $f^n(x_0)$}\label{fig:orbit}
\end{figure}

Under some extra assumptions it is possible to prove that
actually
$$
L^{x_0}(f)\subset F(f).
$$

\section{Groups of isometries}\label{secgp}

\subsection{Generalizations of Hopf's theorem on ends}
The following extends Hopf's theorem that the number of ends of a
finitely generated group is either $0,$ $1,$ $2$, or $\infty$:

\begin{proposition}
\label{prophopf}Assume that $\overline{X}$ is a sequentially compact
Isom($X$)--bordifi\-cation. Let $G$ be a group of isometries fixing a
finite set $F\subset\partial X,$ that is, $GF=F$. If $F$ is not
contained in two stars, then $G$ is bounded.
\end{proposition}

\begin{proof} By passing to a finite index subgroup (which does not
effect the boundedness)\ we can assume that $G$ fixes $F$
pointwise. Now suppose there is a sequence $g_{n}$ in $G$ such
that $g_{n}^{\pm1}x_{0}\rightarrow\xi^{\pm} \in\partial X$. Then
$F$ must be contained in $S(\xi^{+})\cup S(\xi^{-})$ since
otherwise there is is a point in $F$ which on the one hand should
be contracted towards $S(\xi^{+})$ under $g_{n}$, but on the other
hand it is fixed by $G$.
\end{proof}

To see how this implies Hopf's theorem: If two boundary points
belong to different ends, then their stars are disjoint. So if one
has a finitely generated group with finite number of ends, then
applying the proposition with $F$ being the set of ends, one
obtains that the number of ends must be at most two.

By the same method of proof:

\begin{proposition} Assume that $\overline{X}$ is a
sequentially compact Isom($X$)--bordifi\-cation. Let $G$ be a group
of isometries which fixes some collection of stars $S_{i}$ in the
sense that $GS_{i}=S_{i}$ for every $i$. Suppose that for any two
arbitrary stars, there is always an $i$ such that $S_{i}$ is
disjoint from these two stars. Then $G$ is bounded.
\end{proposition}

These two statements can be useful to rule out the existence of
compact quotients of certain Riemannian manifolds or complex
domains.

\subsection{Commuting isometries and free subgroups} The proof of
Proposition \ref{propfixed} in fact shows the following:

\begin{proposition}\label{propcomm} Let $g$ be an isometry and
$\overline{X}$ an Isom($X$)--bordification. Suppose that
$g^{n_i}x_0\rightarrow\xi\in\partial X$ and let $Z(g)$ denote the
centralizer of $g$ in Isom($X$). Then $Z(g)S(\xi)=S(\xi)$,
$Z(g)F(\xi)=F(\xi)$, and $Z(g)F(g)=F(g)$ (when it exists).
\end{proposition}

\begin{proposition}
\label{propfree}Assume that $\overline{X}$ is compact. Let $g$ and $h$ be two
isometries such that $g^{\pm n_{k}}\rightarrow\xi^{\pm}\in\partial
X$, $h^{\pm m_{l}}\rightarrow\eta^{\pm}\in\partial X$ for some
subsequences $n_k$ and $m_l$. Assume that $S(\xi^{+})\cup
S(\xi^{-})$ and $S(\eta^{+})\cup S(\eta^{-})$ are disjoint. Then
the group generated by $g$ and $h$ contains a noncommutative free
subgroup.
\end{proposition}

\begin{proof} By a compactness argument (similar to that in the
proof of Lemma \ref{lemS}) we can find large enough $K$ such that
$$
H(\{g^{n_{k}}x_{0}\}_{k>K})\cup H(\{g^{-n_{k}}x_{0}\}_{k>K})
$$
and
$$
H(\{h^{m_{l}}x_{0}\}_{l>K})\cup H(\{h^{-m_{l}}x_{0}\}_{l>K})
$$
are disjoint. From the contraction observations in subsection
\ref{sseccontr} and the usual freeness criterion \cite{dlH00},
the proposition is proved.
\end{proof}

\begin{figure}[ht!]\small
\psfrag {x+}{$\xi^+$}
\psfraga <-3pt,0pt> {x-}{$\xi^-$}
\psfraga <0pt,3pt> {h}{$\eta^+=\eta^-$}
\cl{\includegraphics{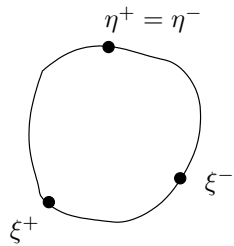}}
\caption{Example of the situation in Proposition \ref{propfree}}
\end{figure}

By a similar proof one has:

\begin{proposition} Assume that $\overline{X}$ is compact. Let
$g$ and $h$ be two isometries such that $g^{\pm
n_{k}}\rightarrow\xi^{\pm}\in\partial X$, $h^{\pm m_{l}
}\rightarrow\eta^{\pm}\in\partial X$ for some subsequences $n_k$
and $m_l$. Assume that $S(\xi^{+}),$ $S(\eta ^{+})\ $and
$S(\xi^{-})\cup S(\eta^{-})$ are disjoint. Then the group
generated by $g$ and $h$ contains a noncommutative free semigroup.
\end{proposition}

\subsection{A metric Furstenberg lemma} The following
can be viewed as an analog of the so-called
Furstenberg's lemma:

\begin{lemma}
\label{lemmetfu}Assume that $\overline{X}$ is a metrizable Isom($X$%
)--compactification
such that $S(\xi)=S^{x_{0}}(\xi)$ for every
$\xi\in\partial X$. Let $g_{n}\in Isom(X)$ and $\mu,\nu$ be two
probability measures on $\partial X.$ Suppose that
$g_{n}\mu\rightarrow\nu$ (in the standard weak topology). Then
either $g_{n}$ is bounded or the support of $\nu$ is contained in
two stars.
\end{lemma}

\begin{proof} We assume that $g_{n}$ is unbounded and by
compactness we select a subsequence so that
$g_{n}^{{}}x_{0}\rightarrow\xi^{+}$, $g_{n}^{-1}x_{0}\rightarrow
\xi^{-}$, and $g_{n}\xi^{-}\rightarrow\xi$. We then have that
$g_{n}S(\xi ^{-})\rightarrow S(\xi)$ in view of the proof of Lemma
\ref{lemscont}. Indeed for any $\eta\in S(\xi^-)$, we have
$g_n\eta\in S(g_n\xi^-)$ and as in the lemma we conclude that
any limit of $g_n\eta$ belongs to $S(\xi)$.

Write $\mu=\mu _{1}+\mu_{2}$ where
$\mu_{1}(\partial X\setminus S(\xi^{-}))=0$ and
$\mu_{2}(S(\xi^{-}))=0$ by letting $\mu_1(A):=\mu(A\cap S(\xi^-))$.
By compactness we can further assume that
$g_{n}\mu_{i}\rightarrow\nu_{i}$ and $\nu=\nu_{1}+\nu_{2}$. Since
$\mu_{1}$ is supported on $S(\xi^{-})$, it follows that $\nu_{1}$
is supported on $S(\xi)$.
Suppose that $f$ is a continuous function vanishing on $S(\xi^{+})$. Then%
$$
\int f(\eta)d\nu_{2}=\lim_{n\rightarrow\infty}\int
f(\eta)d(g_{n}\mu_{2})=\lim_{n\rightarrow\infty}\int f(g_{n}\eta
)d\mu_{2}=0
$$
by the dominated convergence theorem in view of Proposition
\ref{propcontr}. Hence we have shown that supp$\nu\subset
S(\xi)\cup S(\xi^{+})$ as required.
\end{proof}

Furstenberg's lemma, which deals with matrices acting on
projective spaces, has found several beautiful applications since
its first appearance in \cite{Fur63}. For example it is the
key lemma in Furstenberg's
proof of Borel's density theorem, which in turn is
a fundamental tool in the theory of
discrete subgroups in Lie groups.

Our lemma here might be useful for analyzing amenable groups of
isometries (let $\mu=\nu$ be an invariant measure).

\subsection{Random walks}\label{ssecrw}
Let $(X,d)$ be a proper metric space and $\overline{X}$ a
metrizable Isom($X$)--compactifi\-cation. Let $(\Omega,\nu)$ be a
measure space with $\nu(\Omega)=1$ and $L$ a measure preserving
transformation. Given a measurable map $w:\Omega\rightarrow
\rm{Isom}(X)$ we let
$$
u(n,\omega)=w(\omega)w(L\omega)...w(L^{n-1}\omega).
$$
Let $a(n,\omega)=d(x_0,u(n,\omega)x_0)$ and assume that
$$
\int_{\Omega}a(1,\omega)d\nu(\omega)<\infty.
$$
For a fixed $\omega$ we call a subsequence $n_i\rightarrow\infty$
\emph{special} for $\omega$ if there are constants $C$ and $K$ such that
$a(n_i,\omega)>a(m,L^{n_i-m}\omega)-C$ for all $i$ and $m<n_i-K$.
Let $F(\omega)$ denote the face of all limit points of
$u(n,\omega)x_0$ in $\partial X$ along special subsequences.

\begin{theorem}\label{thmrw} Suppose that
\begin{equation}
\liminf_{n\rightarrow\infty}\frac{1}{n}\int_{\Omega}
a(n,\omega)d\nu(\omega)>0.\label{estmulterg}
\end{equation}
Then for a.e.\ $\omega$,
$$
u(n,\omega)x_0\rightarrow \{\eta :F(\omega)\subset S(\eta)\}.
$$
\end{theorem}

\begin{proof} Proposition 4.2 in \cite{KM99} guarantees that
special subsequences exist for a.e.\ $\omega$. From this point on,
the theorem is proved in the same way as
Theorem~\ref{thmindivid}.
\end{proof}

Note that the theorem is not true in general without the
condition (\ref{estmulterg}), since for example
a.e.\ trajectory  of the standard
random walk on the ordinary lattices $\mathbb{Z}^n$ has no
asymptotic direction.

We now specialize to the case when $u(n,\cdot)$ is a random walk
and describe a related result more in terms of boundary theory.

Let $\mathcal{S}$ be the space of closed nonempty subsets of
$\partial X$ with Hausdorff's topology.
Denote by $\Phi(\omega)$ the closure of
$$
\{\xi :\exists C\textrm{ s.t. }
\xi\in\overline{H(u(k,\omega)x_0,C)}\cap\partial X
\textrm{ for all
but finitely many }k \}.
$$
This set may a priori be empty. The next result will guarantee
that it is not empty for a.e.\ $\omega$ and hence we have an
a.e.\ defined map $\omega\rightarrow \mathcal{S}$.
This map is measurable since assigning to a point its
halfspace-closure and the operation of
intersecting closed subsets are continuous.

\begin{theorem}\label{cormu} Let $\mu$ be a probability
measure on a discrete group of isometries $\Gamma$. In the case
$(\Omega,\nu)=\prod_{-\infty}^{\infty}(\Gamma,\mu)$ with $L$ being
the shift, and under assumption (\ref{estmulterg}), the measure
space $(\mathcal{S}$,$\Phi_*(\nu))$ is a $\mu$--boundary of
$\Gamma$.
\end{theorem}

\begin{proof}
Proposition 4.2 in \cite{KM99} guarantees that for a.e.\ $\omega$
there is a $K>0$ and an infinite sequence $n_i$ such that
$$
a(n_i,\omega)>a(n_i-k,L^k\omega)
$$
for all $K<k<n_i$. This means that
\begin{eqnarray*}
d(u(n_i,\omega)x_0,u(k,\omega)x_0)&\leq &
d(u(n_i-k,L^k\omega)x_0,x_0)=a(n_i-k,L^k\omega)\\
&<& a(n_i,\omega)=d(u(n_i,\omega)x_0,x_0)
\end{eqnarray*}
for all $K<k<n_i$ and all $i$. This means that all the limits
of $u(n_i,\omega)x_0$ belong to $\overline{H(u(k,\omega)x_0)}$
for all $k>K$,
in particular $\Phi(\omega)$ is nonempty and belongs to
$\mathcal{S}$.

Consider the path space $\Gamma^{\mathbb{Z}_+}$ with
the induced probability measure $\mathbb{P}$ from the random walk
defined by $\mu$ starting at $e$. Note that $\Gamma$ naturally
acts on $\mathcal F$. The map $\Psi$ gives rise to a map $\Pi$
defined on the path space rather than $\Omega$.

Note that if $\{n_i\}$ is special for $\omega$, then $\{n_{i}-1\}$ is
special for $L\omega$ (cf \cite[page 117]{KM99}).
Special subsequences are moreover
independent of the base point $x_0$. This implies
that $w(\omega)\Psi(L\omega)=\Psi(\omega)$. Now see \cite[1.5]{Kai00}.
\end{proof}

\begin{figure}[ht!]\small
\psfraga <-2pt, 2pt> {x}{$x_0$}
\cl{\includegraphics{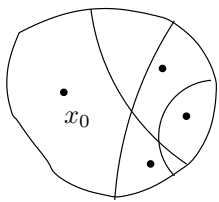}}
\caption{Halfspaces of the random walk intersect.}
\end{figure}

In some situations, eg, when $X$ is $\delta$--hyperbolic
(under some reasonable conditions), the $\mu$--boundary
obtained in the theorem is in fact isomorphic to the Poisson
boundary, see \cite{Kai00}.

The result stated in the introduction follows from the proof of
Theorem \ref{thmrw} and the well-known fact that $A>0$ for simple
random walks on finitely generated nonamenable groups.

\subsection{Proper actions} An isometric action of a
group $\Gamma$ is (\emph{metrically})\emph{ proper} if for every
$x\in X$ and every closed ball $B$ centered at $x$, the set
$\{g\in\Gamma:gx\in B\}$ is finite. A pair of stars $S_1$ and
$S_2$ are \emph{maximal} if the only union of two stars containing
them is $S_1\cup S_2$.

\begin{lemma}\label{thmdiscret} Assume that $\overline{X}$ is
a Hausdorff Isom($X$)--compactification and that $\partial X$ is
not the union of two stars. Suppose that $g$ and $h$ are two
unbounded isometries generating a proper action and that $h^{\pm
n_{j}}x_{0}\rightarrow\xi^{\pm}$ with $S(\xi^{+})$ and $S(\xi
^{-})$ disjoint and maximal. If $g$ fixes $S(\xi^{-})$, then
$h^k=gh^lg^{-1}$ for two nonzero integers $k$ and $l$, and $g$
fixes a star contained in $S(\xi^{+})$.
\end{lemma}

\begin{proof} (Compare \cite{KN02b}.) Since $\overline X$ is a compact
Hausdorff space we can find two disjoint neighborhoods $U^{+}$ and
$U^{-}$ of $S(\xi^{+})$ and $S(\xi^{-})$ respectively, so that
$E:=\overline{X}\setminus(U^{+}\cup U^{-})$ is nonempty and not
contained in $X$. Since $g$ is a homeomorphism fixing
$S(\xi^{\pm})$ we can moreover suppose that
\begin{equation}\label{uinter}
hU^{-}\cap U^{+}=\emptyset.
\end{equation}
Because $h^{-n_{j}}$ contracts toward $S(\xi^{-})$ (Proposition
\ref{propcontrS}) and $g$ is a homeomorphism fixing $S(\xi^{-})$
we have that
$$
gh^{-n_{j}}(E)\subset U^{-}
$$
for all large $j$. In view of (\ref{uinter}) we can find a
$k=k(j)$ such that $h^{k(j)}gh^{-n_{j}}E\cap E$ is nonempty. Let
$g_{j}=h^{k(j)}gh^{-n_{j}}$. Note that
\begin{equation}\label{fixminus}
g_jS(\xi^-)=S(\xi^-)
\end{equation}
and since $g_jS(\xi^+)=h^{k(j)}gS(\xi^+)$, $gS(\xi^+)\cap
gS(\xi^-)=\emptyset$, and $k(j)\rightarrow\infty$,
\begin{equation}\label{fixplus}
g_jS(\xi^+)\rightarrow S(\xi^+).
\end{equation}
In view of (\ref{fixminus}), (\ref{fixplus}), and the assumptions
on $S(\xi^{\pm})$ we have that if $g_{j_k}^{\pm}x_0\rightarrow
\eta^{\pm}\in\partial X$, then either $S(\eta^{\pm})=S(\xi^{\pm})$
or $S(\eta^{\pm})=S(\xi^{\mp})$. In either case this contradicts
that $g_jE\cap E$ is nonempty for all large $j$. Therefore $g_j$
is bounded and by properness we have $g_j=g_i$ for many $i,j$
different. This means that $h^k=gh^lg^{-1}$ for two nonzero
integers $k$ and $l$. Hence
$$
h^kS(g\xi^+)=gh^lg^{-1}gS(\xi^+)=gS(\xi^+)=S(g\xi^+)
$$
and we conclude that $S(g\xi^+)\subset S(\xi^+)$, since
$gS(\xi^-)$ equals all of $S(\xi^-)$.
\end{proof}

\begin{figure}[ht!]\small
\psfrag {h}{$h$}
\psfrag {g}{$g$}
\cl{\includegraphics{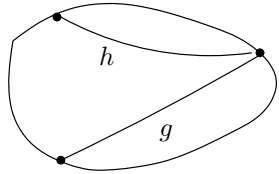}}
\caption{$g$ and $h$ generate a nonproper action.}
\end{figure}

It is instructive to compare Lemma \ref{thmdiscret} with the case
of a Baumslag--Solitar group $<g,h:h^k=gh^lg^{-1}>$ acting on its
Cayley graph.

An \emph{axis} of an isometry is an invariant geodesic line on
which the isometry acts by translation. We say that an isometry
$h$ fixes an endpoint of a geodesic line $c$ if there is a $C>0$
such that $d(h(c(t)),c(t))<C$ for all $t>0$ or all $t<0$.

\begin{proposition}
\label{propaxis}Let $g$, $h$ be two isometries generating a group which acts
properly on a complete metric spaces $X$.
Assume that $g$ has an axis $c$ and that $h$
fixes an endpoint of $c$. Then $[h,g^{N}]=1$ for some $N>0$.
\end{proposition}

\begin{proof} Letting $x_{0}=c(0)$ we have that:
$$
d(x_{0},g^{-n}hg^{n}x_{0})=d(g^{n}x_{0},hg^{n}x_{0})=d(c(nd_{g}),hc(nd_{g}))<C
$$
for all $n>0$ (or $n<0$). As the action of the group is proper, we
must then have that for some $m\neq n$
$$
g^{-m}hg^{m}=g^{-n}hg^{n}
$$
or in other words there is a number $N>0$ such that
$h=g^{-N}hg^{N}$.
\end{proof}

An isometry $g$ is called \emph{strictly hyperbolic} if
$L(g)=\{\xi^{+}\}$ and
$L(g^{-1})=\{\xi^{-}\}$ for two
distinct hyperbolic boundary points $\xi^{+}$ and $\xi^{-}$, which 
by definition means that $S(\xi^{+})=\{\xi^{+}\}$ and
$S(\xi^{-})=\{\xi^{-}\}$.

By the contraction lemma a strictly hyperbolic isometry can have no further
fixed points apart from its two limit points. (Note that if
one knows that $\xi^+$ and $\xi^-$ are hyperbolic limit points,
then it follows that the limit set cannot be larger.)

Examples include pseudo-Anosov elements of mapping class groups,
hyperbolic isometries of a $\delta$--hyperbolic space, and Ballmann's
rank 1
isometries (see \cite{Bal95}, \cite{BB95}) of a CAT(0)--space, see
Proposition \ref{proprank1isom} below.
From Lemma \ref{thmdiscret} and in view of Proposition \ref{propfree}
one has (the star-reflexivity guarantees maximality of any two hyperbolic
boundary points):
\begin{proposition} \label{prop2stricthyp} Assume that $\overline{X}$ is
a Hausdorff star-reflexive Isom($X$)--com\-pact\-ification. The fixed point sets of
two strictly hyperbolic isometries which together generate a proper action
either coincide or are disjoint.
In the latter case, the group generated by the two isometries contains
a noncommutative free subgroup.
\end{proposition}

\section*{Part II\quad Examples and applications}
\addcontentsline{toc}{section}{Part II: Examples and applications}

\section{Hyperbolicity}\label{sechyp}
A boundary point $\xi$ is called \emph{hyperbolic} if $S(\xi)=\{\xi\}$.
A bordification $\overline{X}$ is called \emph{hyperbolic} if
all boundary points are hyperbolic.
A complete metric space $X$ is \emph{asymptotically
hyperbolic} if all stars in $\overline{X}^h$ are disjoint.
It is known that visibility spaces and
Gromov's $\delta$--hyperbolic
spaces (due to P Storm) are asymptotically hyperbolic.

Recall the following standard notation:
$$
(x|z)_{x_0}:=\frac{1}{2}(d(x,x_{0})+d(z,x_{0})-d(x,z)),
$$
and note that
$$
(x|z)\geq\frac{1}{2}d(z,x_{0})
$$
if and only if $x\in H_{z}^{x_0}$, which gives some insight to the
relation between hyperbolicity and halfspaces.

The following axiom is known to hold for the usual boundary of
visibility spaces and
Gromov's $\delta$--hyperbolic spaces, as well as for
the end-compactification, and Floyd's boundary construction
(compare \cite{KN02b}):
\begin{description}
\item[HB] For any $\xi\in\partial X$, there is a family of neighborhoods
$\mathcal{W}$ of $\xi$ in $\overline{X}$, such that the collection of
open sets
$$
\{x:(x|W)>R\}\cup W,
$$
where $W\in\mathcal{W}$, $R>0$, and $(z|W):=\sup_{w\in W\cap X}(z|w)$,
is a fundamental system of neighborhoods of $\xi$ in $\overline{X}$.
\end{description}
\begin{proposition}
\label{prophalfhyp}Every bordification $\overline{X}$ which
satisfies \textbf{HB} is hyperbolic, indeed
$S(\xi)=S^{x_{0}}(\xi)=\{\xi\}$ for every $\xi\in\partial X$.
\end{proposition}

\begin{proof} Given $U$ a neighborhood of $\xi$ in $\overline{X}$
and $C>0$. By definition we may find $R$ and $W\in\mathcal{W}$
such that $\{z:(z|W)>R-C/2\}\subset U$ and by
making $W$ smaller we can also arrange so that $R<d(W,x_{0})/2$
($d(x_{0},\xi)=\infty$). Now
\begin{eqnarray*}
H(W,C)&=&\{z:d(z,W)\leq d(z,x_{0})+C\}\\
&=&\{z:0\leq\sup_{w}\left(  d(z,x_{0})-d(z,w)\right)  +C\}\\
&=&\{z:\inf_{w} d(w,x_0)\leq\sup_{w}\left(  d(z,x_{0})-d(z,w)\right)+
\inf_{w} d(w,x_0)  +C\}\\
&\subset &\{z:d(W,x_{0})\leq\sup_{w}(d(z,x_{0})+d(w,x_{0})-d(z,w))+C\}\\
&=&\{z:(z|W)>R-C/2\}\subset U,
\end{eqnarray*}
which proves the proposition, because $\mathcal{W}$ is a
fundamental system of neighborhoods and $C$ plays no role.
\end{proof}

\begin{figure}[ht!]\small
\psfrag {S}{$S(\xi)=\{\xi\}$}
\cl{\includegraphics{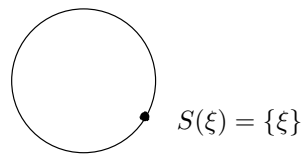}}
\caption{Hyperbolicity}
\end{figure}

For spaces with hyperbolic bordifications, our theory provides
alternative proofs of (mostly) well-known facts, see eg
\cite{BH99} for the theory of ends, \cite{Rat94} for classical
hyperbolic geometry, \cite{Gro87} for word hyperbolic groups, and
\cite{KN02b} for non-locally compact spaces.

\section{Nonpositive curvature}\label{seccat0}

Let $X$ be a complete CAT(0)--space \cite{BH99}.
Recall that the angular metric is
$\angle(\xi,\xi^{\prime})=\sup_{p\in
X}\angle_{p}(\xi^{\prime},\xi)$, where $\xi,\xi^{\prime}$ are
points in the standard visual boundary $\partial X$ of $X$.

\subsection{Stars and Tits geometry}
The
following lemma and its proof can essentially be found in
\cite{BGS 85}:

\begin{lemma}
\label{lemcat0}Let $c$ and $c^{\prime}$ be two geodesic rays emanating from
$x_{0}$ and let $\xi=[c]$ and $\xi^{\prime}=[c^{\prime}]$ be the
corresponding boundary points. Let $p_{i}$ denote the projection
of $c^{\prime}(i)$ onto $c.$ If $\angle(\xi,\xi^{\prime})>\pi/2$
then $p_{i}$ stays bounded as $i\rightarrow\infty$. If
$\angle(\xi,\xi^{\prime})<\pi/2$, then $p_{i}$ is unbounded. In
the case $\angle(\xi,\xi^{\prime})=\pi/2$ then $\{p_i\}$ is
bounded if and only if $x_{0},$ $c,$ and $c^{\prime}$ define a
flat sector.
\end{lemma}

\begin{proof}
First recall the basic angle property of projections 
\cite[Proposition II.2.4]{BH99}: 
$\angle_{p_{i}}(c^{\prime}(i),\xi)\geq\pi/2$ and $\angle_{p_{i}}
(c^{\prime}(i),x_{0})\geq \pi/2$ (when $p_i\neq x_0$). 

If $p_{i}$ is bounded we may assume $p_{i}\rightarrow p$ (along
some subsequence), because the points $p_{i}$ are restricted to a
compact subset of $c$. Then by the upper semicontinuity of angles
(\cite[Proposition II.9.2]{BH99}) we have:
$$
\angle(\xi^{\prime},\xi)\geq\angle_{p}(\xi^{\prime},\xi)\geq\lim\sup_{i}
\angle_{p_{i}}(c^{\prime}(i),\xi)\geq\pi/2.
$$
If $p_{i}$ is unbounded, then in view of \cite[Proposition II.9.8]{BH99} we have%
\begin{eqnarray*}
\angle(\xi,\xi^{\prime})  &
=&\lim_{i\rightarrow\infty}(\pi-\angle_{p_{i}}(c^{\prime}(i),x_{0})-
\angle_{c^{\prime}(i)}(p_{i},x_{0}))\\
&\leq&\pi/2-\lim_{i\rightarrow\infty}\angle_{c^{\prime}(i)}
(p_{i},x_{0})\leq\pi/2.
\end{eqnarray*}
It remains to analyze the case $\angle(\xi^{\prime},\xi)=\pi/2.$
If $p_{i}$ is
a bounded sequence then as above%
\[
\pi/2\geq\angle_{p}(\xi^{\prime},\xi)\geq\lim\sup_{i}\angle_{p_{i}}(c^{\prime
}(i),\xi)\geq\pi/2
\]
and then \cite[Corollary II.9.9]{BH99} shows that $x_{0},$ $c,$ and
$c^{\prime}$ define a flat sector. The converse is trivial:
$p_{i}=x_{0}$.
\end{proof}

\begin{figure}[ht!]\small
\psfrag {x}{$\xi$}
\cl{\includegraphics{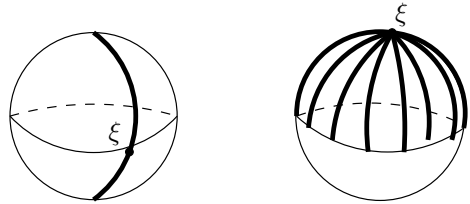}}
\caption{The minimal and maximal star in $\mathbb{R}\times
\mathbb{H}^2$}
\end{figure}

\begin{proposition}
\label{propcat0}Assume $X$ is a complete\/ {\rm CAT(0)}--space and $\overline{X}$ is
the visual bordification. Then
$S(\xi)=S^{x_{0}}(\xi)=\{\eta:\angle (\eta,\xi)\leq\pi/2\}$ for
every $\xi\in\partial X$.
\end{proposition}

\begin{proof} Consider two rays $c_{1}$ and $c_{2}$ from $x_{0}$
representing $\xi$ and $\eta$ respectively. Assume that the
projections of $c_{2}(i)$ onto $c_{1}$ are unbounded. Since by
definition projections realize the shortest distance, we then have
that for any neighborhood $V$ of $\xi$ and for every large enough
$i$ (so that $p_{i}\in V$) that $d(c_{2}(i),V)\leq
d(c_{2}(i),p_{i})\leq d(c_{2}(i),x_{0})$. In the case
$\angle(\eta,\xi)=\pi/2$ and $c_{1}$, $c_{2}$, and $x_{0}$ define
a flat sector, then by Euclidean geometry $V$ contains a point
$\xi^{\prime}$ with $\angle(\xi^{\prime},\eta)<\pi/2$. In view of
Lemma \ref{lemcat0} we hence have
$\{\eta:\angle(\eta,\xi)\leq\pi/2\}\subset S^{x_{0}}(\xi).$

Assume $\angle(\xi,\eta)>\pi/2$ and given $C>0$. By definition
there is a point $y$ such that $\angle_{y}(\xi,\eta)>\pi/2.$ By
continuity (\cite[Proposition II.9.2.(1)]{BH99}) we can find
neigborhoods $V$ of $\xi$ and $U$ of $\eta$ in $\overline{X}$ such
that $\angle_{y}(z,w)\geq\pi/2+\theta$ for every $z\in U,$ $w\in
V$ and some $\theta>0$. Further we make $V$ \ smaller (if
necessary) so that $d(y,V)|\cos(\pi/2+\theta)|\geq
d(x_{0},y)+C^{\prime}$ for some $C^{\prime}>C$. For any $w\in
V\cap X,$ $z\in U\cap X$ we have by the cosine
inequality (ie, comparison with the Euclidean cosine law):%
\begin{eqnarray*}
d(z,w)^{2}  &  \geq &d(y,z)^{2}+d(y,w)^{2}-2d(y,z)d(y,w)\cos\angle_{y}(z,w)\\
&  \geq &d(y,z)^{2}+d(y,w)^{2}+2d(y,z)d(y,w)|\cos(\pi/2+\theta)|\\
&  \geq &
d(y,z)^{2}+(d(x_{0},y)+C^{\prime})^{2}+2d(y,z)(d(x_{0},y)+C^{\prime
})\\
&  =&(d(x_{0},y)+C^{\prime}+d(y,z))^{2}%
\end{eqnarray*}
which implies that $d(z,w)>d(z,x_{0})+C^{\prime}$ by the triangle
inequality. Therefore $d(z,V)>d(z,x_{0})+C$ for all $z\in U\cap X$
and it follows that $\eta\notin S(\xi)$ as desired.
\end{proof}

We will also have use for:

\begin{lemma}\label{lemcat0vis}
Let $X$ be a proper\/ {\rm CAT(0)}--space and assume that $\xi$ is a
hyperbolic point in $\partial X$. Then $\xi$ can be joined to any other
boundary point by a geodesic line in $X$.
\end{lemma}
\begin{proof}
Assume that there is no such geodesic between $\xi$ and $\eta\in\partial
X\setminus \{\xi\}$. By \cite[Theorem 4.11]{Bal95} it then holds that
$\angle (\xi,\eta)\leq\pi$. In fact, there is a (midpoint) $\zeta\in
\partial X$ with $\angle (\xi,\zeta)\leq \pi /2$ (\cite[page 39]{Bal95}),
which contradicts that $S(\xi)=\{\xi\}$ in view of Proposition
\ref{propcat0}.
\end{proof}

\subsection{Corollaries}
All results of the general theory
specialized to the CAT(0)--setting (with the help of Proposition
\ref{propcat0}) seem to be new
except Theorem \ref{thmrw}, Propositions \ref{propcomm} and
\ref{propaxis}. Moreover, in view of Propositions \ref{propcat0}
and \ref{propcontr} (or their proofs in the non-proper case) we
have:

\begin{theorem}\label{thmcat0dyn}
Let $X$ be a complete\/ {\rm CAT(0)}--space. Let $g_{n}$
be a sequence of isometries
such that $g_{n}x_{0}\rightarrow\xi^{+}\in\partial X$ and $g_{n}^{-1}%
x_{0}\rightarrow\xi^{-}\in\partial X$. Then for any
$\eta\in\overline{X}$ with $\angle(\eta,\xi^{-})>\pi/2$ we have
that
$$
g_{n}\eta\rightarrow\{\zeta:\angle(\xi^{+},\zeta)\leq\pi/2\}
$$
(in the sense that $\limsup\angle(\xi^{+},g_n\eta)\leq\pi/2$ when
$X$ is not proper). Assuming that $X$ is proper, the convergence
is uniform outside neighborhoods of $S(\xi^{-})$.
\end{theorem}

Applied to the special case of iterates of a single isometry
$g_n:=h^{k_n}$, the theorem partially extends (since it also deals
with parabolic isometries) a lemma of Schroeder \cite{BGS 85}
generalized by Ruane \cite{Rua01} to include also singular
CAT(0)--spaces. Let us emphasize that this theorem gives
information also about the dynamics of parabolic isometries of
general CAT(0)--spaces.

Combining Propositions \ref{propcat0}\ and \ref{propfree} yields
the following result which generalizes the main theorem in
\cite{Rua01} (because no group is here assumed to act cocompactly
and properly):

\begin{theorem} \label{thmcat0free}
Let $X$ be a proper\/ {\rm CAT(0)}--space. If $g$ and
$h$ are two unbounded isometries with limit points $\xi^{-}$,
$\xi^{+}$ and $\eta^{-}$, $\eta^{+}$ respectively (not necessarily
all distinct), with $Td(\{\xi^{\pm}\},\{\eta^{\pm}\})>\pi$, then
the group generated by $g$ and $h$ contains a noncommutative free
subgroup.
\end{theorem}

The following proposition generalizes \cite[Lemma 4.5]{BB95} and
\cite[Theorem 8]{Swe99} (by weakening the hypothesis):

\begin{proposition} Let $X$ be a complete\/ {\rm CAT(0)}--space and $g$
a hyperbolic isometry with an axis $c$. Assume that $h$ is an
isometry which fixes one endpoint of $c$ and that $g$ and $h$
generate a group acting properly. Then $h$ commutes with some power
of $g$ and $h$ fixes both endpoints of $c$.
\end{proposition}

\proof First note that the notion of fixing an
endpoint of a geodesic coincide with the usual one for CAT(0)--spaces
and the standard ray boundary $\partial X$.  From Proposition
\ref{propaxis} we have $h=g^{-N}hg^{N}$. Therefore
$$
h(c(\pm\infty))=\lim_{n\rightarrow\infty}hg^{\pm
nN}x_{0}=\lim_{n\rightarrow\infty}g^{\pm nN}hx_{0}=c(\pm\infty).\eqno{\qed}
$$

An isometry is called a rank 1 isometry if it is hyperbolic
with an axis which does not bound a flat halfplane \cite{Bal95}.
The usefulness of this notion was demonstrated by Ballmann and
collaborators.
\begin{proposition}\label{proprank1isom}
Rank 1 isometries are strictly hyperbolic.
\end{proposition}
\begin{proof}
Recall that $g$ fixes the stars of its limit points. So unless
$S(\xi^{\pm})=\{\xi^{\pm}\}$ this would contradict the contraction lemma
for rank 1 isometries \cite[Lemma III.3.3]{Bal95}.
\end{proof}

Actually, the converse is also true since strictly hyperbolic
isometries clearly cannot be elliptic, and also not parabolic
(look at preserved horoballs) and that the axis cannot bound
a flat halfplane in view of Proposition \ref{propcat0}.
We obtain the following theorem which sheds some light on the
question \cite[Question III.1.1]{Bal95}, see also \cite[Theorem
III.3.5]{Bal95}:

\begin{theorem} \label{thmballfree}
Suppose that a group $\Gamma$ acts properly
by isometry on a proper\/ {\rm CAT(0)}--space. If the limit set
contains at least three points, one of which is hyperbolic,
then $\Gamma$ contains noncommutative free subgroups.
\end{theorem}
\begin{proof}
Let $\xi$ be a hyperbolic boundary point. Then together with
any other boundary point it does not bound a flat halfplane
in view of Lemma \ref{lemcat0vis} and Proposition \ref{propcat0}.
Therefore given a sequence $g_n$ in $\Gamma$ which we can assume
that $g_nx_0\rightarrow \xi$ and $g_n^{-1} x_0\rightarrow \eta$
for some other boundary point $eta$, which we moreover suppose is
different from $\xi$. Indeed, if $\xi =\eta$ then by the basic
contraction lemma and since the limit set contains at least three
points, we can find such $g_n$.

The lemma \cite[Lemma III.3.2]{Bal95} now guarantees the existence
of a rank 1 isometry $g$ with hyperbolic limitpoints say $\xi^{\pm}$.
By assumption there is another
limit point $\eta$ for the group different from
$\xi^{\pm}$, take $h_n$ for which $h_nx_0\rightarrow \eta$. Since
the boundary is star-reflexive, by the basic contraction lemma
we have that some $h_N$ moves one of $\xi^{\pm}$ say $\xi^{+}$.
Moreover it moves a neighborhood of $\xi^{+}$ into a neighborhood
of the star of $\eta$. Now consider the sequence $h_Ng^n$ (in
$n$), this has one limit point outside both $\xi^{\pm}$ and the
other, the hyperbolic point $\xi^{-}$. Apply again Ballmann's
lemma to obtain another rank 1 isometry $h$. Since it does not
have have $\xi^{+}$ as limit point, the theorem is proved in view
of Proposition \ref{prop2stricthyp}.
\end{proof}

In the end it might be more powerful to use ping-pong
arguments with the halfspaces directly without pushing it to
the boundary. For example,
one can in this way extend the main theorem in
\cite{AFN02}
somewhat: the condition of no-fake-angles can be removed and the
translation lengths do not necessarily have to be \emph{strictly}
greater than the length of $S$.

\section{Hilbert's geometry on convex sets}

Let $X$ be a bounded convex domain in $\mathbb{R}^{n}$ and
$\partial X$ the usual boundary. The Hilbert metric on this 
domain is a complete metric and is defined as follows. For any 
two distinct points $x$
and $y$ draw the chord through these points. Now $d(x,y)$
is the logarithm of the projective cross-ratio of $x$, $y$, and 
the two endpoints of the chord. 
We refer to \cite{dlH93} or
\cite{Nuss88} for more information, note in particular that
semicontractions of Hilbert's metric arise in several situations,
for example in potential theory. Recall that in this context the
\emph{star} of a boundary point $\xi$, Star($\xi$), is the
intersection of $\partial X$ with the union of all hyperplanes
which are disjoint from $X$ but contain $\xi$. We have:

\begin{proposition}\label{prophilb} Assume that $X$ is a bounded convex domain
equipped with Hilbert's metric and let $\overline{X}$ be the
closure in $\mathbb{R}^n$. Then
$S(\xi)=S^{x_{0}}(\xi)=$Star($\xi)$ for every $\xi\in\partial X.$
\end{proposition}

\begin{proof}
The inclusion $S(\xi)\subset$Star$(\xi)$ follows from the inclusion%
\[
H(W,C)\subset\{z:(z|W)\geq\frac{1}{2}d(W,x_{0})+C^{\prime}\}
\]
proved in Proposition \ref{prophalfhyp} using the same
terminology, together with the proof of Theorem 5.2 in
\cite{KN02}. The other inclusion follows because given $W$ and
$\zeta$ it is simple to see that we can approximate $\zeta$ with a
point arbitrary far from $x_{0}$ but staying on finite Hilbert
distance to $W$ (the Hilbert metric remains bounded near a line
segment of the boundary in the direction parallel to this line
segment). In particular, $\bigcap\overline{H(V,C)}$ is independent
of $C$ and equals $S^{x_0}(\xi)$.
\end{proof}

\begin{figure}[ht!]\small
\psfrag {h}{$\eta$}
\psfrag {z}{$\zeta$}
\psfrag {x}{$\xi$}
\cl{\includegraphics{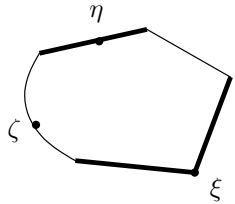}}
\caption{Examples of stars for Hilbert's metric}
\end{figure}

For these metric spaces, it seems that several results obtained
in this paper cannot be found in the literature. For example, we
obtain from combining Propositions \ref{prophopf} and \ref{prophilb}:
\begin{theorem}\label{thmhilb}
Any polyhedral cone with noncompact automorphism group has
simplicial diameter at most 3.
\end{theorem}.

\begin{figure}[ht!]
\cl{\includegraphics{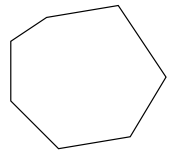}}
\caption{A convex set with compact automorphism group}
\end{figure}

The simplicial diameter is the smallest number of simplices required
to connect any two points.
In dimension $2$ a rather complete result concerning which convex sets
have nonfinite automorphism group can be found in \cite{dlH93}.


The literature on symmetric
or homogeneous cones is vast. For recent works on cones where the
automorphism group admits a cocompact lattices, see the works
of Y Benoist. eg \cite{Ben01}.
Hilbert's metric can also be a tool in the study of Coxeter groups
via the Tits cone such as in \cite{McM01}.

\section{Several complex variables}

Let $X$ be a bounded domain in $\mathbb{C}^{N}$. We denote by $d_{X}$ the
Kobayashi distance and by $F_{X}$ the corresponding infinitesimal metric on
$X$. The metric space $(X,d_{X})$ is not always complete (pseudoconvexity $X$
is for example a necessary condition), but when it is, $(X,d_{X})$ is in
addition proper and geodesic. We refer to \cite{Kob98} for more details.

\subsection{Points of strict pseudoconvexity}
Here is a relation between strictly pseudoconvex points and
stars:
\begin{theorem}\label{thmkob} Let $X$ be a bounded domain in
$\mathbb{C}^n$ with $C^2$--smooth boundary equipped with
Kobayashi's metric. If $\xi_1$ and $\xi_2$ are two distinct
boundary points at which $X$ is strictly pseudoconvex, then
$s(\xi_1,\xi_2)\geq 2$.
\end{theorem}
\begin{proof}Combining \cite[Theorem 4.5.8]{Kob98} with an estimate
due to Forstneric--Rosay, cf \cite[Corollary 4.5.12]{Kob98}, one
has for some constant $C$ and fixed $x_0$,
$$
d(z_1,z_2)\geq C+d(z_1,x_0)+d(z_2,x_0)
$$
for all $z_1$ (resp.\ $z_2$) sufficiently close to $\xi_1$ (resp.\
$\xi_2)$. Hence $\xi_1\notin S(\xi_2)$ and $\xi_2\notin S(\xi_1)$.
\end{proof}

\begin{corollary}Let $X$ be a strictly pseudoconvex bounded
domain with $C^2$--boundary equipped with Kobayashi's metric. Then
$S(\xi)=S^{x_0}(\xi)=\{\xi\}$ for every $\xi\in\partial X$.
\end{corollary}

\subsection{From metric to distance estimates}\label{ssecmetdist}

Let $X$ be a bounded $C^{1+\alpha}$--smooth ($\alpha>0$) domain in
$\mathbb{C}^{N}$ which is complete in the Kobayashi metric and fix some
$x_{0}\in X$.  Euclidean distances are denoted by $\delta$.
Assume that for some $\varepsilon>0$ and $c_{1}>0$%
\begin{equation}
F_{X}(z;v)\geq c_{1}\frac{||v||}{|\delta(z,\partial X)|^{\varepsilon}}
\label{est1}%
\end{equation}
for all $z\in X$ and $v\in\mathbb{C}^{N}$. Examples include bounded
pseudoconvex domains with real analytic boundary \cite{DF 79}
and $C^{2}$--strictly
pseudoconvex domains, see Theorem E.3, and section X.10.4 in \cite{JP 93}.

\begin{lemma}\label{lemvis}
Let $\gamma$ be a minimizing geodesic between two points $z$ and $w$ in $X$.
Then%
\begin{equation}
\delta(z,w)\leq C(2d_{X}(x_{0},\gamma)+\varepsilon^{-1})e^{-\varepsilon
d_{X}(x_{0},\gamma)} \tag{$*$}\label{estlemvis}%
\end{equation}
for some $C>0$ depending only on $X$ and $x_{0}$.
\end{lemma}

\begin{proof}
Let $m$ be a point on $\gamma$ of minimal distance $r:=d_{X}(x_{0},\gamma)$ to
$x_{0}$ and denote by $\gamma_{1}:[0,a]\rightarrow X$ the (reparametrized)
piece of $\gamma$ going from $m$ to $z$. Because of the minimality of $r$ and
the triangle inequality we have%
\begin{align*}
d_{X}(x_{0},\gamma_{1}(t))  &  \geq r\\
d_{X}(x_{0},\gamma_{1}(t))  &  \geq t-r
\end{align*}
for all $t$. The following estimate is known (in the case of
$C^{2}$--smoothness see \cite[Theorem 4.5.8]{Kob98} or
\cite[X.10.4]{JP 93}, and in the more general case it is a
consequence of \cite[Proposition 2.5]{FR 87}): there is a constant
$c_{3}$ such that
\begin{equation}
d_{X}(x_{0},z)\leq c_{3}-\log\delta(z,\partial X) \label{est2}%
\end{equation}
for all $z\in Z$.

In the case $a>2r$, we have from the above estimates, since $\gamma_{1}$ is a
unit speed geodesic that:%
\begin{align*}
\delta(m,z)  &  \leq\int_{0}^{a}||\dot{\gamma}_{1}(t)||dt\leq c_{1}\int
_{0}^{a}\delta(\gamma_{1}(t),\partial X)^{\varepsilon}F_{X}(\gamma_{1}%
(t);\dot{\gamma}_{1}(t))dt\\
&  =c_{1}\int_{0}^{a}\delta(\gamma_{1}(t),\partial X)^{\varepsilon}dt\leq
c_{4}\int_{0}^{a}e^{-\varepsilon d_{X}(x_{0},\gamma_{1}(t))}dt\\
&  \leq c_{4}\int_{0}^{2r}e^{-\varepsilon r}dt+c_{4}\int_{2r}^{a}%
e^{-\varepsilon(t-r)}dt<c_{4}2re^{-\varepsilon r}+c_{4}\varepsilon
^{-1}e^{-\varepsilon r}.
\end{align*}
In the case, $a\leq2r$ we make the same estimate but without decomposing the
integral. By a symmetric argument with $w$ instead of $z$, the lemma is proved
in view of the triangle inequality.
\end{proof}

\begin{theorem}
\label{thmfloyd}The closure $\overline{X}$ is a hyperbolic compactification of
$X$, indeed $S(\xi)=S^{x_0}(\xi)=\{\xi\}$ for every $\xi\in\partial X.$
\end{theorem}

\begin{proof}
First note that the right hand side of (\ref{estlemvis}) in Lemma
\ref{lemvis} tends to $0$ if
$d_{X}(x_{0},\gamma)\rightarrow\infty$ ($x_{0}$ is fixed). Now recall the
simple and standard fact that (for any geodesic space)%
\[
(z_{1}|z_{2})_{x_{0}}:=\frac{1}{2}(d_{X}(z_{1},x_{0})+d_{X}(z_{2},x_{0}%
)-d_{X}(z_{1},z_{2}))\leq d_{X}(x_{0},\gamma)
\]
for any geodesic segment joining $z_{1}$ and $z_{2},$ see eg \cite{KN02}.
This means that the condition \textbf{HB} and the assertion follows
from Proposition \ref{prophalfhyp}.
\end{proof}

We record the following result which is formulated in a more traditional style
but which we have not been able to find in the literature.

\begin{theorem}
Given
two distinct boundary points $\xi_{1}$, $\xi_{2}\in\partial X$ there exists
constants $\kappa>0$ and $c\in\mathbb{R}$ depending only on $X$, $\xi_{1}$ and
$\xi_{2}$ such that%
\[
d_{X}(z_{1},z_{2})\geq c+\log\frac{1}{\delta(z_{1},\partial X)}+\log\frac
{1}{\delta(z_{2},\partial X)}.
\]
\end{theorem}

\begin{proof}
In view of the proof of Theorem \ref{thmfloyd} and Lemma \ref{lemvis}
we have that for some
neighborhoods of $\xi_{1}$ and $\xi_{2}$ there is a constant $R$ such that for
any $z_{1}$ and $z_{2}$ in these neighborhoods respectively,
$$
(z_{1}|z_{2})_{x_{0}}\leq2R
$$
which spelled out reads%
\begin{align*}
d(z_{1},z_{2})  &  \geq R-d(z_{1},x_{0})-d(z_{2},x_{0})\\
&  \geq c+\log\frac{1}{\delta(z_{1},\partial X)}+\log\frac{1}{\delta
(z_{2},\partial X)}%
\end{align*}
in view of (\ref{est2}).
\end{proof}

\subsection{Convex domains}

Recall that the face $\mathrm{Face}(\xi)$ is the intersection of
all hyperplanes which contains $\xi$ but avoids the interior of
the convex set. The following result is due to Abate (see \cite{A
90} or \cite[Corollary 2.4.25]{A 89}):

\begin{theorem}
Let $X$ be a convex $C^{2}$--smooth bounded domain in $\mathbb{C}^{N}$ and
given $\xi_{1}$, $\xi_{2}\in\partial X$ such that $\xi_{1}\notin F(\xi_{2})$
(and hence also $\xi_{2}\notin F(\xi_{1})$). Then there exists $\kappa>0$ and
$c\in\mathbb{R}$ such that%
\[
d_{X}(z_{1},z_{2})\geq c+\log\frac{1}{\delta(z_{1},\partial X)}+\log\frac
{1}{\delta(z_{2},\partial X)}%
\]
for any $z_{1}\in X\cap\{w:\delta(w,F(\xi_{1}))<\kappa\}$ and $z_{2}\in
X\cap\{w:\delta(w,F(\xi_{2}))<\kappa\}$.
\end{theorem}

\begin{corollary}
Let $X$ be a convex $C^{2}$--smooth bounded domain in $\mathbb{C}^{N}$. Then
$S(\xi)\subset \mathrm{Star}(\xi)$ for any $\xi\in\partial X$.
\end{corollary}

\begin{proof}
This is deduced similarily to Theorem \ref{thmkob}.
\end{proof}

What are the stars for a general
bounded pseudoconvex (Kobayashi hyperbolic) domain with Kobayashi's metric?
Note here that
Hilbert's metric is an analogous metric and Teichm\"uller metric
is another example.

\subsection{Iteration of holomorphic maps on bounded domains}

The iteration of holomorphic self-maps of bounded domains in $\mathbb{C}$ was
studied by Wolff, Denjoy, Valiron and Heins. In several variables perhaps the
first works were done by H Cartan and Herv\'{e}. From 1980 and
onwards there have appeared many papers by several authors including
Vesentini, Abate, Vigue, D Ma, X-J Huang, Zhang, Ren, and Mellon.
Most often the main tool is an appropriately
generalized Schwarz--Pick lemma.

Since holomorphic maps semicontract Kobayashi distances,
$$
d(f(x),f(y))\leq d(x,y)
$$
for all $x,y\in X$, we can canonically associate the set $F(f)$
to any holomorphic map
(in view of Theorem \ref{thmindivid}). From Theorem \ref{thmkob},
we obtain the following:

\begin{theorem}Let $X$ be a $C^2$ bounded domain in
$\mathbb{C}^n$, $f:X\rightarrow X$ a holomorphic map, and $d$ the
Kobayashi distance. Assume that $(X,d)$ is complete.
Then $F(f)$ contains at most one point of
strong pseudoconvexity.
\end{theorem}

\begin{figure}[ht!]\small
\psfraga <0pt,-2pt> {x}{$x_0$}
\psfraga <-5pt,0pt> {F}{$F(f)=\xi$}
\cl{\includegraphics{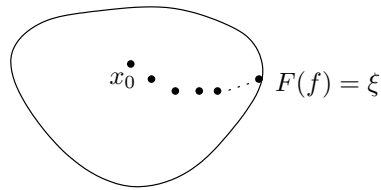}} \caption{Wolff--Denjoy
result: orbit of $f$}
\end{figure}

To see the relevance of this to iterations of holomorphic maps,
recall that (Theorem \ref{thmindivid}) if the stars are only singletons, then
for any $z$
either $f^m(z)$ stays away from $\partial X$ ($F(f)=\emptyset$),
or $\lim_{m\rightarrow\infty} f^m(z)=\xi$ for some $\xi\in\partial
X$ ($F(f)=\{\xi\}$).
In particular, we can formulate the following new Wolff--Denjoy
result (the case of
$n=2$ and real analytic boundary, was proved by Zhang and Ren
in \cite{ZR95}):

\begin{theorem}\label{thmWD}
Let $X$ be a bounded domain in $\mathbb{C}^n$ such as in
subsection \ref{ssecmetdist}. Let $f:X\rightarrow X$
be a holomorphic map. Then either the orbit of $f$ stays away from the
boundary or there is a unique boundary point $\xi$ such that
$$
\lim_{m\rightarrow\infty} f^m(z)=\xi
$$
for any $z\in X$.
\end{theorem}

Finally we remark that that in view of \cite{Ka01} it would be
interesting to identify the horofunction boundary for domains with
Kobayashi's metric.

\section{Teichm\"uller spaces and mapping class groups}

Let $M$ be a closed surface
of genus $g\geq2$ and denote by $\mathcal{T}$ the associated
Teichm\"{u}ller space $\mathcal{T}$. It is known that $\mathcal{T}$
can be embedded as a bounded domain in
$\mathbb{C}^{N}$ and that it is pseudoconvex but not strictly
pseudoconvex. By a theorem of Royden, one also knows that Kobayashi's metric
coincides with Teichm\"{u}ller's metric.

\subsection{Stars in the Thurston boundary}\label{ssecteich}
Let $\mathcal{S}$ be the set of homotopy
classes of simple closed curves on $M$. Denote by $i(\alpha
,\beta)$ the minimal number of intersection of representatives of
$\alpha,\beta\in\mathcal{S}.$ Let $\mathcal{MF}$ (resp.\
$\mathcal{PMF}$) be the set of (resp.\ projective equivalence
classes of) measured foliations, which coincides with the closure
of the image of the embedding
$$
\alpha\mapsto i(\alpha,\cdot)
$$
of $\mathcal{S}$ into $\mathbb{R}^{\mathcal{S}}$ (resp.\
$P\mathbb{R}^{\mathcal{S}}$). The intersection number $i$ extends
to a bihomogeneous continuous function on
$\mathcal{MF}\times\mathcal{MF}$. A foliation $F\in\mathcal{PMF}$
is called \emph{minimal} if $i(F,\alpha)>0$ for every
$\alpha\in\mathcal{S}$. The Teichm\"{u}ller space $\mathcal{T}$ of
$M$ is embedded into $P\mathbb{R}^{\mathcal{S}}$ by the hyperbolic
length function. (Below, $V_{\phi}$ stands for the vertical
foliation and $h_{\phi}$ the horizontal length associated to a
quadratic differential $\phi$, see \cite{KM96} for more details.)

\begin{lemma}
\label{lemteich}Let $\phi_{n}$ be the quadratic differential corresponding (in
the Teichm\"{u}ller embedding with reference point $x_{0}$) to
$x_{n}\in \mathcal{T}$. Assume $\phi_{n}\rightarrow\phi_{\infty}$,
a norm one quadratic differential, and $x_{n}\rightarrow F$ in
$\mathcal{PMF}$. Whenever $\beta _{n}\in\mathcal{S}$ such that
$Ext_{x_{n}}(\beta_{n})<D$ and $\beta _{n}\rightarrow H$ in
$\mathcal{PMF}$, it holds that
$$
i(V_{\phi_{\infty}},H)=0=i(F,H).
$$
\end{lemma}

\begin{proof} A proof analysis shows that this is proved in
\cite{Mas82}: Denote by $\psi_{n}$ the terminal differential
corresponding to the Teichm\"{u}ller map from $x_{0}$ to $x_{n}$.
Since
$$
h_{\psi_{n}}(\beta_{n})\leq Ext_{x_{n}}(\beta_{n})^{1/2}<D^{1/2},
$$
$x_{n}\rightarrow\infty$ in $\mathcal{T}$, and in view of the
stretching of the Teichm\"{u}ller map
($h_{\psi_{n}}=e^{d(x_{0},x_{n})}h_{\phi_{n}}$) we see that
$$
\lim_{n\rightarrow\infty}i(V_{\phi_{n}},\beta_{n})=\lim_{n\rightarrow\infty
}h_{\phi_{n}}(\beta_{n})=0.
$$
Since $\beta_{n}$ is a sequence in $\mathcal{S}$ converging to $H$
in $\mathcal{PMF}$, there is a sequence $\lambda_{n}$ of bounded
positive scalars such that $\lambda_{n}\beta_{n}\rightarrow H$ in
$\mathcal{MF}$. By continuity and homogeneity of $i$ we have
$i(V_{\phi_{\infty}},H)=0.$

For the second equality note that it is known that
$x_{n}\rightarrow F$ in $\mathcal{PMF}$ implies that there is a
sequence $r_{n}\rightarrow0$ such that
$i(r_{n}x_{n},\cdot)\rightarrow i(F,\cdot)$ in $\mathcal{MF}$.
From definitions we also have
$$
i(x_{n},\beta_{n})\leq
A_{g}^{1/2}Ext_{x_{n}}(\beta_{n})^{1/2}<A_{g}^{1/2}D^{1/2},
$$
which by the same argument as before now also shows the second
equality.
\end{proof}

The following result can be viewed as a slight generalization of Lemma
1.4.2 in \cite{KM96} although formulated in a very different way 
(note that there seems that there is a misprint 
in their statement) and is obtained by almost the same proof.

\begin{theorem} \label{thmteich}
Let $X$ be the Teichm\"{u}ller space of a
compact surface and equipped with the Teichm\"{u}ller metric $d$.
Let $\overline{X}$ be the Thurston compactification
$X\cup\mathcal{PMF}$. For $F\in\mathcal{PMF}$, a minimal
foliation, we have
$$
S(F)\subset\{G:i(F,G)=0\}.
$$
\end{theorem}

\begin{proof} Given $y_n\rightarrow G\in S^{x_0}(F)$, select
$x_{n}\rightarrow F$ such that $d(y_{n},x_{n})\leq
d(y_{n},x_{0})+C$ for all $n$ and some $C>0$. From continuity and
Mumford compactness, it is a fact that sequences $\beta_{n}$ as in
Lemma \ref{lemteich} corresponding to $x_n$ always exist. Assume
now that $F$ is minimal. It is then known (due to Rees) that,
$i(F,G)=0$ if and only if $G$ is minimal and equivalent to $F$.
Hence $V_{\phi_{0}}$, $F$ and $H$ as in Lemma \ref{lemteich} are
all equivalent minimal foliations. Fix these. Note that
$\lambda_{n}\rightarrow0$ here because of the minimality. Let
$\theta_n$ (resp.\ $\psi_n$) denote the initial (resp.\ terminal)
quadratic differential of the Teichm\"uller map from $x_0$ to
$y_n$. We have
\begin{eqnarray*}
i(V_{\theta_{n}},\lambda_{n}\beta_{n})  &
=&\lambda_{n}h_{\theta_{n}}(\beta_{n})\\
&  =&\lambda_{n}e^{-d(y_{n},x_{0})}h_{\psi_{n}}(\beta_{n})\\
& \leq
&\lambda_{n}e^{-d(y_{n},x_{0})}De^{d(y_{n},x_{n})}\rightarrow 0,
\end{eqnarray*}
where the last inequality follows from Kerckhoff's formula for
Teichm\"{u}ller distances. Thus $i(V_{\theta_{0}},H)=0$, which
implies what we want, since $i(F,G)=0$ is an equivalence relation
for minimal foliations and because of Lemma \ref{lemteich}.
Finally since the set on the right in the theorem is closed,
we have $i(F,G)=0$ for all $G\in S(F)$.
\end{proof}

The set of uniquely ergodic foliations is denoted by $\mathcal{UE}$
and is a subset of full Lebesgue measure in the Thurston boundary.

\begin{corollary}\label{corue}
Every point $F\in\mathcal{UE}$ is a hyperbolic boundary point, indeed,
$S(F)=S^{x_{0}}(F)=\{F\}$.
\end{corollary}

\begin{conjecture}\label{conjteich}
For any $F\in\mathcal{PMF}$, it holds that
$$
S(F)=\{G:i(F,G)=0\}.
$$
\end{conjecture}

The conjecture would imply that the Thurston boundary equipped with
the star distance restricted to $\mathcal{S}$ is the 1--skeleton of
the curve complex (it is interesting to here recall the important
result of Masur--Minsky that this complex is Gromov hyperbolic).

\subsection{Mapping class groups}
Although the arguments
in this paper provide (especially if all stars of the
Teichm\"uller spaces can be identified) an alternative explanation
of some theorems on the mapping class groups of surfaces obtained
notably in \cite{Iva87} and \cite{MP89}, it might however be
preferable to study the action directly on the Thurston boundary
(or the curve complex) as is done in those works.

It is a standard fact that the only fixed points of a pseudo-Anosov
element are two uniquely ergodic foliations and so from Corollary \ref{corue}
we get:
\begin{proposition}
Pseudo-Anosov elements of the mapping class groups
are strictly hyperbolic.
\end{proposition}

Hence Proposition \ref{prop2stricthyp} applies and gives known facts.
Theorem \ref{thmteich} and the fundamental contraction
property gives a new, more analytic approach, as well as
some additional information,
to Theorems 7.3.A and 7.3.B in \cite{Iva01}:

\begin{theorem}\label{thmteichdyn}
Suppose that $g_n$ is a sequence of elements in the mapping class groups
for which $g_n^{\pm 1}x_0$ converge to two minimal foliations $F^{\pm}$
in $\overline{X}$,
the Thurston compactification. Then for $z$ outside a neighborhood of
$\{G:i(F^-,G)=0\}$ in $\overline{X}$, $g_nz$ converges uniformly to
$\{G:i(F^+,G)\}$.
\end{theorem}

\section{Infinite groups}

This section contains some remarks and questions concerning
infinite groups from the point of view of stars and dynamics
of isometries.

\textbf{Free subgroups}\qua 
Let $\Gamma$ be a group with a left invariant word metric $||\cdot||$.
So $||g||=d(g,e)$ where $e$ is the identity element in $\Gamma$. Moreover,
$||xy^{\pm 1}||=\min\{d(x,y),d(x,y^{-1}\}$. We may formulate the following
freeness criterion:

\begin{lemma}
Let $g$ and
$h$ be two elements of order at least $3$ in $\Gamma$ and 
let $\Lambda$ be the subgroup
generated by $g$ and $h$. If for any $a\in\Lambda$ at least one of
$||ag^{\pm 1}||$ and $||ah^{\pm 1}||$ is strictly greater than
$||a||$, then $\Lambda\cong F_{2}$.
\end{lemma}

\begin{proof}
The statement $||ag||>||a||$ is equivalent to that $a\notin H^{e}_{g^{-1}}$.
So  $||ag^{\pm 1}||>||a||$ means that $a\notin  H^{e}_{g}\cup H^{e}_{g^{-1}}$.
Therefore the hypothesis implies that $H^{e}_{g}\cup H^{e}_{g^{-1}}$ and
$H^{e}_{h}\cup H^{e}_{h^{-1}}$ are disjoint inside the invariant set $\Lambda$. 
In view of the first observation in
subsection \ref{sseccontr} and the ping-pong lemma
\cite{dlH00} applied to the unions of the halfspaces 
associated to $g^{\pm 1}$ and $h^{\pm 1}$ respectively, and intersected
by $\Lambda $,
the lemma now follows.
\end{proof}

There is a similar criterion
for free semigroups.

\textbf{Random walks}\qua 
Let $\Gamma$ be a finitely generated group, $\mu$ the uniformly
distributed probability measure on a finite generating set $A$,
$X$ the Cayley graph associated to $A$ and $\overline X^h$ the
horofunction compactification. If $\Gamma$ is nonamenable, then
Theorem \ref{cormu} provides a (probably often nontrivial)
$\mu$--boundary for $(\Gamma,\mu)$. In particular we have that
if the the random walk has a linear rate of escape, then, from
some time on, all the
halfspaces defined by the points of the random walk intersect.
This is a nontrivial phenomenon, which
for example can be seen from thinking
about random walks on $\mathbb{Z}^n$.

\textbf{Associated incidence geometries}\qua
Let $\Gamma$ denote a finitely generated group with a boundary
$\partial \Gamma$.  Consider the incidence geometry (cf
\cite{Tit74}) defined by the points and the stars in
$\partial\Gamma$, and acted upon by $\Gamma$. In general or for
some specific group, what can this geometry be? In view of section
3, for example torsion, subexponential growth, or amenability
of $\Gamma$ implies strong restrictions.

\begin{figure}[ht!]
\cl{\includegraphics{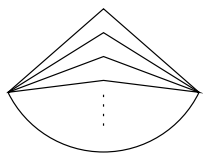}} \caption{Incidence
geometry of the Baumslag--Solitar group BS(1,$n$)}
\end{figure}

As we have seen above, examples of metric spaces and associated
incidence geometry are:
\begin{itemize}
\item Gromov hyperbolic spaces -- trivial
\item Hilbert's geometry -- the face lattice
\item CAT(0)--space -- Tits geometry
\item mapping class groups -- the curve complex (conjectural).
\end{itemize}

It may also be interesting to
extend the existing theory of convergence groups to a more general
setting where one has nontrivial incidence geometry mixed in. A
first instance of what we essentially have in mind can be found in
\cite{FM99}.

\textbf{Rigidity theory}\qua 
The following philosophy, vaguely formulated, lies behind the Mostow--Margulis
rigidity theory. Any proper homomorphism of one group to another
gives rise to an incidence preserving map at infinity of the groups.
Incidence preserving maps at infinity must be of a very special
kind. Can one make this more precise and when does it (or a part of it)
hold? The most simple cases are homomorphisms from $\mathbb{Z}^n$ into the
isometry group of some hyperbolic space. Another instance (but now
with trivial incidence geometry) is the Floyd--Cannon--Thurston maps
in Kleinian group theory.

\textbf{Group cohomology}\qua
The following is relevant for (the first) $\mathrm{L}^2$--coho\-mology:
Consider the class of harmonic (ie, satisfying the mean-value property)
functions on the vertices of an oriented
Cayley graph such
that its differential (which is a function on the edges, the difference
of the values on the two vertices) is square summable. These functions are
called Dirichlet harmonic functions. Compactify the
graph in the Stone--Cech way relative to this family of functions.
The group does
not admit nonconstant Dirichlet harmonic functions if and only if
this compactification is the one-point compactification (an easy
consequence of the maximum principle). What is the star geometry of
this compactification? For many groups it seems that the compactification
should be hyperbolic.

\end{document}